\documentclass[12pt]{amsart}
\usepackage[backref=page,colorlinks,bookmarksopen=true]{hyperref}
\usepackage{pdfsync}
\usepackage{amssymb}
\usepackage{latexsym}
\usepackage{amsmath}

\def\surl#1_#2{\mathrel{\mathop{\kern 0pt #1}\limits_{#2}}}

\newcommand{\eenmatrix}[4]{\setlength{\unitlength}{1ex}
\raisebox{-1.5ex}[0ex][3.5ex]{\parbox{7.5ex}{\begin{picture}(5,5)(0,0)
\put(1.8,1.9){\line(0,1){4}}
\put(1.8,5.9){\line(1,0){4}}
\put(5.8,5.9){\line(0,-1){4}}
\put(5.8,1.9){\line(-1,0){4}}
\put(2.3,6.2){\parbox[b]{3ex}{\centering $#1$}}
\put(6.2,3.4){\parbox[b]{25ex}{$#2$ }}
\put(-3,3.4){\parbox[b]{1.5ex}{  \raggedleft $#3$  }}
\put(-1,0){\parbox[t]{10ex}{\centering $#4$}}
\end{picture}

}}}

\begin{document}

\title{ Measured quantum groupoids associated with matched pairs of locally compact groupoids}
\author{Jean-Michel Vallin}
\address{Administrative address:  MAPMO  UMR 6628 CNRS Universit\'{e} d'Orl\'eans, France  \  \  \   \  \  \  \  \  \ Alternative address: IMJ \\Plateau 7D, 175
rue du Chevaleret 75013 Paris, France  
\\ Email: \href{mailto:jmva@math.jussieu.fr}{\textcolor{blue}{jmva@math.jussieu.fr}}  }


\markboth{Jean-Michel Vallin} {Matched pairs of locally compact groupoids}

\begin{abstract} Generalizing the notion of matched pair of groups, we define   and study matched pairs of locally compact groupoids endowed with Haar systems, in order to give new examples of measured quantum groupoids. \end{abstract}
\maketitle
\tableofcontents
 \newpage

\section{Introduction}
\def\I{ I_{\mathcal H,\mathcal K}}
\def\k{\gamma^K}
\def\h{\delta^H}
\newenvironment{dm}{\hspace*{0,15in} {\it Proof:}}{\qed}

 Dealing with locally compact groupoids,  we already defined in  the articles \cite{Val0} and \cite{Val4}) a notion of pseudo multiplicative unitary and Hopf bimodule in order to generalize, in that framework, classical notions of multiplicative unitary  (\cite{BS}) and Hopf von Neumann algebras (\cite{ES}) which led to locally compact quantum groups (\cite{BS}, \cite{W}, \cite{KV}......).
\vskip 0.2cm
In an other article (\cite{EV}), starting with any  depth 2
inclusion of von Neumann algebras $M_0 \subset  M_1$, with an operator-valued weight $T_1$ verifying
a regularity condition,  Michel Enock and the author  have constructed a pseudo-multiplicative unitary generating two Hopf bimodules in duality; one of them acts on $M_1$ in such a way  that $M_0$ is isomorphic to the fixed point algebra and  the von Neumann algebra $M_2$, given by the basic construction,  is isomorphic to the crossed product.
\vskip 0.2cm

 The axiomatic of locally compact quantum groupoids has been developed by Franck Lesieur in \cite{L1}  and \cite{L2} and simplified by M.Enock (\cite{E2} Appendice),  who has also studied  the theory of their actions on von Neumann  algebras, generalizing previous results due to  S.Vaes (\cite{Vas2}).
\vskip 0.2cm
The aim of this  article is to give a large number of examples of measured quantum groupoids as defined by M.Enock and F.Lesieur. We generalize   at the same time the notion of matched pair of finite groupoids (\cite{Val2}) and of locally compact  group (\cite{BS}, \cite{VV},  \cite{BSV}, \cite{BSV1}....), in order  to obtain  such examples  coming from a suitable pseudo multiplicative unitary. 

 In the second paragraph are  recalled general definitions  about locally compact quantum groupoids and their actions on   von Neumann algebras.  

We specify, in the third chapter,  the notion of  matched pair of locally compact groupoids, we prove that, for such a pair, 
there exists a canonical  action of each groupoid on the other one,   and   we give families of examples. Finally  we find  a canonical pseudo multiplicative unitary which generates their crossed products.
  
  In the fourth chapter we investigate the Hopf bimodule structures of the  crossed products given by the pseudo multiplicative unitary, and find suitable Haar operator valued weights for these    structures. 
  
  We  study,  in the last chapter, two kinds of examples. The first one is pretty natural and comes from  matched pairs of groups   actions:   a very general example of matched pair of groups is  the "ax+b" group  (\cite{BSV1} Chap 4), and  pentagonal transformations   lead also to such actions ( \cite{BSV1} Prop 5.1).  We prove that, for any    matched pair $G_1G_2$ of locally  compact groups  in the sense of \cite{BSV},   which acts on a locally compact space $X$,   then $X\times G_1$, $X\times G_2$  is a matched pair of groupoids. Moreover,    $G_1$ (resp. $G_2$)   acts, as a group, on  the space $X\times G_2$ (resp.$X\times G_1$)  in such a way that their    usual crossed product   is   the one obtained  using chapter 3. We investigate  the quantum groupoid structure  given to these crossed products by chapter 4 which is actually different from the one given to any crossed product as   dual of a transformation group. The second example is the farthest possible from groups, it      comes from principal groupoids of the form $X\times X$ where $X = X_1 \times X_2$ is the cartesian product of two locally compact spaces $X_1$ and $X_2$, we prove that the structures given by the previous chapters mixes the ones  given by the pair groupoids $X_1 \times X_1$ and $X_2 \times X_2$.
\vskip 0.2cm
Several continuations of this article can  be considered. One can  weaken the condition $\mathcal G_1 \cap \mathcal G_2 = \mathcal G^0$,  which even with  finite groups or groupoids  gives substantial examples (\cite{Val3},  \cite{AN} 2.8). Also  
   a characterization of these objects  in terms of cleft extensions in the spirit of S.Vaes and L.Vainerman \cite{VV} should be obtained.

\section{ Measured quantum groupoids and their actions}

\subsection{ Measured  quantum groupoids }
\label{intro}

 Let's recall the definition of a measured quantum groupoid  due to  M.Enock which   extends  F.Lesieur's works.  We use \cite{L1}, \cite{L2} and \cite{E2} for  general references, in particular we suppose known spatial theory and relative tensor products (\cite{C}, \cite{S}) .
 
\subsubsection{\bf{Definition}}
\label{mqg}
{\it  A {\bf measured quantum groupoid} is     a special collection   

$\mathfrak G = (N,M,\alpha,\beta,\Gamma,T,T',\nu)$  such that:

i) $M, N$ are two  von Neumann algebras,  $\alpha: N \to M$ and $\beta:  N^o \to M$ are  commuting  faithful normal non degenerate representations,
\vskip0.3cm
 ii)  $\Gamma: M \to  M\surl{\nonumber_{\beta} \star_{\alpha}}_{N}M $ is a one to one normal morphism such that: 

 \centerline{$\Gamma(\beta(x))=1 \surl{\nonumber_{\beta} \otimes_{\alpha}}_{N}
  \beta(x)$ }
  
 \centerline{$\Gamma(\alpha(x))= \alpha(x) \surl{\ _{\beta}
    \otimes_{\alpha}}_{N} 1 $}

  \centerline{$(\Gamma \surl{\nonumber _{\beta}
    \star_{\alpha}}_{N} id)\Gamma=(id \surl{\ _{\beta} \star_{\alpha}}_{N} \Gamma)\Gamma.$}

\vskip0.3cm
iii)  $T $  (resp $T'$) is   a faithful semi finite normal operator valued weight  from $M$ to $\alpha(N)$ (resp $\beta(N)$) such that: 
\vskip0.2cm
 \centerline{ $(id \surl{\nonumber_{\beta} \star_{\alpha}}_{N} T )\Gamma(x) =  T(x) \surl{\nonumber_{\beta} \otimes_{\alpha}}_{N} 1$ for any $x \in \mathcal M^+_{T}$}

 \centerline{$(T'\surl{\nonumber_{\beta} \star_{\alpha}}_{N} i)\Gamma(x) =  1 \surl{\nonumber_{\beta} \otimes_{\alpha}}_{N} T'(x)$ for any $x \in \mathcal M^+_{T'}$}

\vskip0.3cm
iv) $\nu$ is a faithful semi finite normal weight on $N$ which is relatively invariant  with respect  to $T $ and $ T'$, i.e.   for any $t \in \mathbb R:  \sigma^{\Phi }_t\sigma^{\Psi}_t =  \sigma^{\Psi}_t\sigma^{\Phi }_t $, where $\Phi = \nu \circ \alpha^{-1} \circ T$ and $\Psi = \nu \circ \beta^{-1} \circ T'$ .

}

\subsubsection{\bf{Remark}}
\label{vacances}
The assertion iii) can be replaced by the weights conditions:

iii)'   $(id \surl{\nonumber_{\beta} \star_{\alpha}}_{N} \Phi )\Gamma(x) =  T(x)  $ for any $x \in \mathcal M^+_{T}$ 

\hskip 0.8cm$(\Psi \surl{\nonumber_{\beta} \star_{\alpha}}_{N} i)\Gamma(x) =   T'(x)$ for any $x \in \mathcal M^+_{T'}$

We shall say that the quantum groupoid is commutative (resp.symmetric) if $M$ is abelian (resp. $\varsigma \Gamma = \Gamma $, where $\varsigma : M\surl{\nonumber_{\beta} \star_{\alpha}}_{N}M  \to M\surl{\nonumber_{\alpha} \star_{\beta}}_{N^o}M$ is the natural  flip).
\subsubsection{\bf{Definition (\cite{EV} 5.6)}}
\label{upm}
{\it  Let $N$ be a von Neumann algebra and $\nu$ a faithful normal semifinite weight on $N$, let $\alpha$ (resp. $\beta, \hat {\beta}$) be a faithful non degenerate representation (resp.   anti representations) on a Hilbert space $\mathfrak H$ commuting two by two, a {\bf pseudo multiplicative unitary} $W$ over the basis {\bf $(N,\alpha,\beta,\hat \beta)$} is  a unitary from $ \mathfrak H \surl{\nonumber_{{\beta}} \otimes_{\alpha}}_{\nu} \mathfrak H$ to $ \mathfrak H \surl{\nonumber_{\alpha} \otimes_{\hat \beta}}_{\nu^o} \mathfrak H$ such that: 
\begin{itemize}
\vskip0.3cm
\begin{item} $W$ intertwines $\alpha,\beta,\hat \beta$,  which means that for any $n \in N$ one has: 
\vskip0.1cm
 \centerline{$W(\alpha(n) \surl{\nonumber_{ \beta} \otimes_{\alpha}}_{N} 1) = ( 1 \surl{\nonumber_{\alpha} \otimes_{\hat \beta }}_{N^o} \alpha(n)) W $ }
  
 \centerline{$W(1 \surl{\nonumber_{  \beta} \otimes_{\alpha}}_{N}   \beta (n)) = ( 1 \surl{\nonumber_{\alpha} \otimes_{ \hat \beta }}_{N^o}  \beta(n)) W $}

  \centerline{$W(\hat \beta(n) \surl{\nonumber_{  \beta} \otimes_{\alpha}}_{N} 1) = (\hat  \beta(n) \surl{\nonumber_{\alpha} \otimes_{ \hat \beta }}_{N^o} 1) W $}
  
  \centerline{$W(1 \surl{\nonumber_{  \beta} \otimes_{\alpha}}_{N} \hat \beta(n)) = (  \beta(n) \surl{\nonumber_{\alpha} \otimes_{\hat \beta }}_{N^o} 1)) W $}

\end{item}
\vskip0.3cm
\begin{item} The operator $W$ satisfies the "pentagonal" equation: 
\vskip0.2cm
\begin{multline*}
(1_{\mathfrak H}  \surl{\nonumber_{\alpha} \otimes_{\hat \beta}}_{N^o} W )(W  \surl{\nonumber_{ \beta} \otimes_{\alpha}}_{N} 1_{\mathfrak H}) = \\
= (W  \surl{\nonumber_{ \beta} \otimes_{\alpha}}_{N^o} 1)(\sigma_{\nu^o}  \surl{\nonumber_{\alpha} \otimes_{\hat \beta}}_{N^o} 1_{\mathfrak H})(1_{\mathfrak H}  \surl{\nonumber_{\alpha} \otimes_{\hat \beta}}_{N^o} W)\sigma_{2\nu}(1_{\mathfrak H } \surl{\nonumber_{  \beta} \otimes_{\alpha}}_{N} \sigma_{\nu^o})(1_{\mathfrak H}  \surl{\nonumber_{  \beta} \otimes_{\alpha}}_{N} W)
\end{multline*}
where $\sigma_{\nu^o} $  is the flip map: $\mathfrak H  \surl{\nonumber_{\alpha} \otimes_{\hat \beta}}_{\nu^o} \mathfrak H \to \mathfrak H  \surl{\nonumber_{\hat \beta} \otimes_{\alpha}}_{\nu} \mathfrak H $ and $\sigma_{2\nu}$ is the flip map:  $\mathfrak H  \surl{\nonumber_{\alpha} \otimes_{\hat \beta}}_{\nu^o} \mathfrak H \surl{\nonumber_{\alpha} \otimes_{\hat \beta}}_{\nu^o} \mathfrak H \to \mathfrak H  \surl{\nonumber_{\hat \beta} \otimes_{\alpha}}_{\nu} \mathfrak H \surl{\nonumber_{\alpha} \otimes_{\hat \beta}}_{\nu^o} \mathfrak H $
\end{item}
\end{itemize}}

\subsubsection{\bf{Remark}} 
\label{echocard}
In fact measured quantum groupoids and  pseudo multiplicative unitaries are  closely linked. According to    \cite{EV} chap. 6,  if $W$ is a pseudo multiplicative unitary on $\mathcal L(\mathfrak H)$,   it generates two von Neumann algebras $M$ (its right leg) and $\widehat M$ (its left leg) and two coproducts $\Gamma$ and $\widehat \Gamma$ on $M$ and $\widehat M$ respectively,  i.e. two maps verifying definition \ref{mqg} ii). More precisely, for any $m \in M$ and $\hat m \in \widehat M$, one has: $\Gamma(m) = W^*(1 \surl{\nonumber_{  \alpha} \otimes_{\hat \beta}}_{N^o} m)W$ and $\widehat \Gamma( \hat m) = \sigma_{\nu^o}W(\hat m \surl{\nonumber_{  \beta} \otimes_{\alpha}}_{N}1)W^*\sigma_{\nu}$;  and conversely for a given measured quantum groupoid, one can associate a pseudo multiplicative unitary  to it,  with a manageability condition (implying weak regularity) which leads to a duality theory and the following theorem: 

\subsubsection{\bf{Theorem}}(\cite{L1},\cite{L2},\cite{E2})
\label{lesieur}
{\it Let $\mathfrak G =   (N,M,\alpha,\beta,\Gamma,T ,T',\nu)$ be a measured quantum groupoid,   for any  $n$ in $N$ let's define $\hat \beta(n)= J_\Phi \alpha(n^*)J_\Phi$,  then one can associate  to $\mathfrak G$ a pseudo multiplicative unitary $W$, over the basis  $(N, \alpha,\beta, \hat \beta) $,  independent of the choice of $T$ and $T'$,  the left  leg of which is $M$ and gives  $\Gamma$.  The dual coproduct $\hat \Gamma$ on the  right leg $\widehat M$ leads to a  measured quantum groupoid $(N, \hat M, \alpha, \hat \beta, \widehat T, \widehat T',\nu)$ denoted     $\widehat{\mathfrak G}$, the dual of $\mathfrak G$.}

To the quantum groupoid $\mathfrak G$ is associated a canonical $*$-antiautomorphism of $M$ called the coinverse, which is involutive, and  verifies the condition $R \circ \alpha = \beta$ and $\Gamma \circ R = \varsigma_{N^o}(R\surl{\nonumber_{\beta} \otimes_{\alpha}}_{N} R)\Gamma$.

Let $j$  be the application  defined for any $x \in M$ by $j(x) = J_\Phi x^*J_\Phi$,  and set $\Gamma^c = (j \surl{\nonumber_{ \beta} \star_{\alpha}}_{N}j)\Gamma \circ j,T^c =j \circ T\circ j ,R^c = j\circ R\circ j$.

Then one can consider      two other  quantum groupoids:
\begin{itemize}
\begin{item}
the commutant  $\mathfrak G^c = (N^o,M',\hat \beta,\hat \alpha,\Gamma^c,T^c,R^cT^cR^c,\nu^0)$ 
\end{item}
\begin{item}
 the commutant of the dual $ (\widehat{\mathfrak G})^c = (N^o,(\widehat M)', \beta,\hat \alpha,(\widehat \Gamma)^c,(\widehat T)^c,(\widehat R)^c(\widehat T)^c(\widehat R)^c,\nu^0)$ 
this last is  an important tool for the    duality of actions.
\end{item}
\end{itemize}

\vskip 0,4cm

\subsection{Measured quantum groupoids in action}

As quantum groups act on von Neumann algebras, measured quantum groupoids act on (von Neumann ) modules with isomorphic  basis. Generalizing in this  context what we have done in the finite dimensional situation in \cite{Val2}, M.Enock has given in \cite{E2} a nice framework for these actions together with   double crossed product theorems. Let's recall some of his definitions.

\subsubsection{\bf{Definition}}
\label{debut}
{\it Let $\mathfrak G =  (N,M,\alpha,\beta,\Gamma,T_L,T_R,\nu)$ be a measured quantum groupoid, and let $A$ be a von Neumann algebra acting on a hilbert space $\mathcal H$. A right (resp. left) action of $\mathfrak G$ on $A$ is a pair $(b,\mathfrak a)$ such that:

i) $b: N \to A$ is an injective $*$-antihomomorphism  (resp.morphism),

ii) $\mathfrak a: A \to A  \surl{\nonumber_{b} \star_{\alpha}}_{N} M$ (resp. $ A  \surl{\nonumber_{b} \star_{\beta}}_{N^o} M$)  is an injective $*$-homomorphism,

iii)  for all $n \in N$ one has: $\mathfrak a(b(n)) = 1  \surl{\nonumber_{b} \otimes_{\alpha}}_{N} \beta(n)$ (resp.$\mathfrak a(b(n)) = 1  \surl{\nonumber_{b} \otimes_{\beta}}_{N^o} \alpha(n)$)   and one has:  $(\mathfrak a \surl{\nonumber_{b} \star_{\alpha}}_{N} id)\mathfrak a = ( 1  \surl{\nonumber_{b} \star_{\alpha}}_{N} \Gamma)\mathfrak a $ (resp. $(\mathfrak a \surl{\nonumber_{b} \star_{\beta}}_{N^o} id)\mathfrak a = ( 1  \surl{\nonumber_{b} \star_{\beta}}_{N^o} \varsigma \Gamma)\mathfrak a $).}

\subsubsection{\bf{Definitions}}
\label{codebut}
{\it  Let $(b,\mathfrak a)$ be a right (resp.left) action of a given measured quantum groupoid $\mathfrak G$ on a given von Neumann algebra $A$, then:

i) the crossed product of $A$  by the action $(b,\mathfrak a)$ is  the   sub von Neumann algebra of $ A  \surl{\nonumber_{b} \star_{\alpha}}_{N} \mathcal L(\mathcal H_\phi)$ (resp.$ A  \surl{\nonumber_{b} \star_{\beta}}_{N^o} \mathcal L(\mathcal H_\phi)$) generated by  $\mathfrak a(A)$ and $1   \surl{\nonumber_{b} \otimes_{\alpha}}_{N} \widehat M'$ (resp. $1   \surl{\nonumber_{b} \otimes_{\beta}}_{N^o} \widehat M$), it will be denoted by $A \rtimes_{\mathfrak a}\mathfrak G$.  

ii) the invariant subalgebra  is defined by:
$$A^{\mathfrak a} = \{  x \in A \cap b(N)' / \mathfrak a(x) = x \nonumber_{b} \otimes_{\alpha} 1\}$$
 }

\vskip 0,4cm

\subsubsection{\bf{Theorem}(\cite{E2} 6.13, 9.3, 9.5, 10.11, 10.12)}
\label{stefaan}
{\it  Let $(b,\mathfrak a)$ be a right   action of a given measured quantum groupoid $\mathfrak G$ on a given von Neumann algebra $A$,  and set   $\Phi = \nu \circ \alpha^{-1}\circ T$, then:

i) for any $x \in A^+$,  the extended positive element of $A$:
$$ T_{\mathfrak a}(x) =  (id  \surl{\nonumber_{b} \star_{\alpha}}_{\nu} \Phi)\mathfrak a(x) $$
is an extended positive element of $A^{\mathfrak a}$ and $T_{\mathfrak a}$ is a normal faithful operator valued weight from $A$ to  $A^{\mathfrak a}$

ii) there exists a unique action $( 1\nonumber_{b} \otimes _{\alpha} \hat \alpha,\tilde {\mathfrak a} ) $ of $\ \hat{\mathfrak G}^c$ on $A \rtimes_{\mathfrak a} \mathfrak G$ which verifies  for any $x \in A$, $y \in \hat M'$:
$$ \tilde {\mathfrak a}(\mathfrak a(x)) = \mathfrak a(x)\nonumber_{\hat \alpha} \otimes _{\beta} 1$$
$$ \tilde {\mathfrak a}(1 \nonumber_{b} \otimes _{\alpha}y ) = 1 \nonumber_{b} \otimes _{\alpha}\hat \Gamma^c(y)$$

iii) for any $y \in \hat M'$: $ T_{\tilde{\mathfrak a}}(  1 \nonumber_{b} \otimes _{\alpha}y  ) = 1 \nonumber_{b} \otimes _{\alpha}\hat T^c(y) =   \mathfrak a(b(\beta^{-1}(\hat T^c(y))))$, $ T_{\tilde{\mathfrak a}}$ is semi finite, and  $(A \rtimes_\alpha \mathfrak G)^{\tilde{\mathfrak a}} = \mathfrak a(A)$.}

\subsubsection{\bf{Corollary and definition}}
\label{poidsdual}
{\it For any normal semi finite faithful operator valued weight  $\theta$ on $A$,  the operator valued weight $\tilde{\theta} = \mathfrak a \circ \theta \circ \mathfrak a^{-1} \circ T_{\tilde{\mathfrak a}}$ is normal semi finite faithful on $A \rtimes_{\mathfrak a}\mathfrak G$ and  will be called the {\bf dual operator valued weight} of $\theta$.}

\vskip 0,7cm
In fact, the  examples we shall deal with in this article, come from an action of a commutative  measured quantum groupoid on  a commutative von Neumann algebra, nevertheless, as we shall see, the crossed product will be a substantial non commutative quantum groupoid. So let us  see more precisely  the commutative situation.

\subsection{The abelian case} 
\label{R}
 All  measured quantum groupoids, involved in this article, will have a commutative basis, and the Hilbert spaces will be of the form $L^2(X, dx)$ where $X$ is a second countable locally compact space endowed with a Radon measure, hence this will simplify the  relative tensor products. Let  $(Y,dy)$ be a locally compact space  endowed with a Radon measure and let $(L^2(X,dx),\beta)$  (resp. $(L^2(Z,dz), \alpha)$) be a faithful normal representation of $N= L^\infty(Y,ds)$, then the relative tensor product   $  L^2(X,dx) \surl{\nonumber_{\beta} \otimes_{\alpha}}_{N} L^2(Z,dz)$ is the completion of   the algebraic tensor product  $\mathcal K(X)  \odot \mathcal K(Z)$ equipped with the pre-scalar product, defined for any  $f_1,f_2 \in \mathcal K(X)$ and any $g_1,g_2 \in \mathcal K(Z)$ by:
 \begin{align*}
  (f_1 \odot g_1, f_2 \odot g_2) 
  &= (\alpha(\frac{d\omega_{f_1,f_2}\circ \beta}{dy})g_1,g_2) = (\beta(\frac{d\omega_{g_1,g_2}\circ \alpha}{dy})f_1,f_2) \\
  &= \int_{Z} \alpha(\frac{d\omega_{f_1,f_2}\circ \beta}{dy})g_1(z)\overline{g_2(z)}dz = \int_{X} \beta(\frac{d\omega_{g_1,g_2}\circ \alpha}{dy})f_1(x)\overline{f_2(x)}dx
  \end{align*}
In our framework, the relative tensor product $  L^2(X,dx) \surl{\nonumber_{\beta} \otimes_{\alpha}}_{N} L^2(Z,dz)$  will  also be  viewed as $L^2(X_\beta \times_\alpha Z, dx  _\beta{\times_\alpha} dz)$,  where $X_\beta \times_\alpha Z$  is  a suitable fibred product of $X$ and $Z$ and $dx  _\beta{\times_\alpha}dz$ is a Radon measure; we shall intensively use this identification throughout  this paper.

The commutative measured quantum groupoids (i.e $M$ is commutative) are, as expected, coming from measured groupoids in the sense of  Jean Renault (\cite{R},\cite{AR}). Notwithstanding the fact that a Weyl theorem does not exist  in that context, up to some inessential reduction (see \cite{Ra} theorem 4.1)   we can deal with a Hausdorff  locally compact groupoid $\mathcal G$, we shall suppose it is $\sigma$- compact and endowed  with a Haar system $\{ \lambda^u / u \in \mathcal G^0\}$ and a quasi invariant measure $\nu$ on $\mathcal G^0$, we shall denote $\mu = \int_{\mathcal G^0} \lambda^ud\nu(u)$ the integrated measure on $\mathcal G$, for any $u \in \mathcal G^0$, $\lambda_u$ will be the image of  $\lambda^u$ by the application $g \mapsto g^{-1}$, and $\delta$ will be the Radon Nikodym derivative of $\mu^{-1}$ w.r.t.  $\mu$ . This will allow us to use the $*$-algebra of continuous numerical functions on $\mathcal G$ with compact support, which will be noted $\mathcal K(\mathcal G)$. 

If $s(g) = g^{-1}g$ (resp. $r(g) = gg^{-1}$) is the source   (resp. goal) of any element $g$ of $\mathcal G$, then  one can define two representations of $N= L^\infty(\mathcal G^0,\nu)$ on $M = L^\infty(\mathcal G,\mu)$, (which are also antirepresentations)   defined for any $f \in L^\infty(\mathcal G^0,\nu)$ by:
$$s_{\mathcal G}(f) = f \circ s, \  \  r_{\mathcal G}(f) = f \circ r$$

One easily verifies that for any $f,f' \in \mathcal K(\mathcal G)$ and for $\nu$-almost any $u \in \mathcal G^0$, one has:
$$ \frac{d\omega_{f ,f'}\circ r}{d\nu}(u) = \int_{\mathcal G }f(g) \overline{f' (g)}d\lambda^u  \  \ , \  \  \frac{d\omega_{f ,f'}\circ s}{d\nu}(u) = \int_{\mathcal G}\delta(v)f(v^{-1})\overline{f'(v^{-1})}d\lambda^u(v)$$

\subsubsection{\bf{Remark and notations}}
\label{gilardi}
One can define for any $i \in \{s ,r \}$,   $\mathcal G^2_{i,j} = \{(g,g') \in \mathcal G \times \mathcal G / i(g) = j(g') \}$,  for instance  $\mathcal G^2_{s, r } = \mathcal G^{2}$ (the set of pairs of composable elements), let us note  $\mu^2_{i,j} = \mu_{i_{\mathcal G}} \times _{j_{\mathcal G}}\mu$,  so  $L^2(\mathcal G,\mu) \surl{\nonumber_{i_{\mathcal G}} \otimes_{j_{\mathcal G}}}_{L^\infty(\mathcal G^0,\nu)}L^2(\mathcal G,\mu)$ and $L^2(\mathcal G^2_{i,j},\mu^2_{i,j})$ are  isomorphic and isometric. 
\vskip 0.5cm

In \cite{Val0} and \cite{Val4}  we proved that we  can associate to  $\mathcal G$ a   pseudo multiplicative unitary $W_{\mathcal G}: L^2(\mathcal G^2_{s ,r },\mu^2_{s ,r }) \to L^2(\mathcal G^2_{r ,r },\mu^2_{r ,r })$. 
It is given for any $\xi \in L^2(\mathcal G^2_{s ,r },\mu_{s ,r })$ and $\mu^2_{r  ,r }$ almost any  $(x,y) \in G^2_{r ,r }$ by:
$$W_{\mathcal G}\xi(x,y) = \xi(x,x^{-1}y)$$ 

The left leg of $W_{\mathcal G}$ generates the  commutative quantum groupoid:
$$\mathfrak G(\mathcal G) = (L^\infty(\mathcal G^0,\nu), L^\infty(\mathcal G,\mu), r_{\mathcal G},s_{\mathcal G}, \Gamma_{\mathcal G}, T_{\mathcal G},  T^{-1}_{\mathcal G} ,\nu),$$ with the following formulas, for any $f \in \mathcal K(\mathcal G)$:

\begin{itemize}

\item \   \    \   \    \    \  {\bf Coproduct}

$\Gamma_{\mathcal G}(f)(x,y) = f(xy)\ \  \mathrm{for}   \  \  \mathrm{ any}   \ \  (x,y) \in G^2_{s_{\mathcal G},r_{\mathcal G}} $

\item \   \    \   \    \    \  {\bf Left and right operator valued weights}

For $\mu$-almost any element $g \in \mathcal G$:
\begin{align*}
T_{\mathcal G}(f)(g) 
&= \int_{  \mathcal G} f(x)d\lambda^{r(g)}(x)\\
  T^{-1}_{\mathcal G}(f)(g) &= \int_{  \mathcal G} f(y)d\lambda_{s(g)}(y)
\end{align*}
\end{itemize}

The right  leg of $W_{\mathcal G}$ generates the symmetric quantum groupoid, which means the coproduct is invariant by a natural flip.
$$\widehat {\mathfrak G}(\mathcal G)= (L^\infty(\mathcal G^0,\nu), \mathcal L(\mathcal G), r_{\mathcal G}, r_{\mathcal G},\widehat \Gamma_{\mathcal G}, \widehat T_{\mathcal G}, \widehat T_{\mathcal G},\nu),$$
 where $\mathcal L(\mathcal G)$ is the left regular algebra of $\mathcal G$, which is the sub von Neumann algebra of $\mathcal L(L^2(\mathcal G,\mu))$ generated  by the operators defined for any $f,h\in \mathcal K(\mathcal G)$ by:
$$ \lambda(f)h(x) = \int_{\mathcal G} f(g)h(g^{-1}x)d\lambda^{r_{\mathcal G}(x)}(g)$$
For any $f \in \mathcal K(\mathcal G)$, $\xi \in L^2(G^2_{r_{\mathcal G}, r_{\mathcal G}},\mu_{r_{\mathcal G}, r_{\mathcal G}})$ and almost any $(x,y) \in \mathcal G^2_{r_{\mathcal G}, r_{\mathcal G}}$, one has:

\begin{itemize}

\item \   \    \   \    \    \  {\bf Coproduct}:

$\hat \Gamma_{\mathcal G}(\lambda(f))\xi(x,y) =  \int_{\mathcal G} f(g)\xi(g^{-1}x,g^{-1}y)d\lambda^{r_{\mathcal G}(x)}(g) $

\item \   \    \   \    \    \  {\bf Left  (=  right ) operator valued weight}

$\widehat T_{\mathcal G}(\lambda(f)) =  r_{\mathcal G}(f\mid \mathcal G^0)$

\end{itemize}

\subsubsection{\bf{Remark}}
\label{porcine}
Of course one can consider the {\bf right} regular representation of $\mathcal G$ which generates in $\mathcal L(L^2(\mathcal G,\mu))$ the commutant of $\mathcal L(\mathcal G)$ and gives a commutant structure of quantum groupoid: 
$$\widehat {\mathfrak G}'(\mathcal G)= (L^\infty(\mathcal G^0,\nu), \mathcal R(\mathcal G), s_{\mathcal G}, s_{\mathcal G}, \widehat \Gamma_{\mathcal G}', \widehat T_{\mathcal G}', \widehat T_{\mathcal G}',\nu),$$

\subsubsection{\bf{The  pair groupoid example}} 
\label{fima1}
 We suppose that $\mathcal G = X \times X$  where $X$ is an Hausdorff  locally compact space together with some Radon  measure $\nu$. $\mathcal G$ is given its natural locally compact groupoid structure, $\mathcal G^0$ is the diagonal of $X\times X$ which will be identified with $X$, its Haar system  is $(\delta_x\otimes \nu)_{x \in X}$ and $\nu$ is a (quasi) invariant measure.

The von Neuman algebra $\mathcal R(\mathcal G)$ is isomorphic to $\mathcal L(L^ 2(X,\nu))$ as for any $h\in \mathcal K(X^2)$, one has $\rho(h) = 1 \otimes T_h$, where $T_h$ is the integral operator defined for any $ \xi \in L ^2(X, \nu)$ by: $T_h\xi(y) = \int_X h(y,b)\xi(b)d\nu(b)$

Here the space $L^2(\mathcal G^2_{s,r},\mu^2_{s,r})$ (resp.$L^2(\mathcal G^2_{s,s},\mu^2_{s,s})$) can be identified with $L^2(X^3,\nu^{\otimes 3})$ using the application: $((x,y),(y,t)) \to (x,y,t)$ from $\mathcal G^2_{s,r}$ to $X^3$ (resp.$((x,y),(t,y)) \to (x,y,t)$ from $\mathcal G^2_{s,s}$ to $X^3$). With this identification, for any $f \in \mathcal K(\mathcal G)$, and any $(x,y) \in \mathcal G$, one has:
$$   \Gamma_{\mathcal G}(f) = f_{13} , \hskip 1.5cm \widehat \Gamma_{\mathcal G}'(\rho(f)) = 1 \otimes T_f\otimes 1$$
$$   T_{\mathcal G}(f))(x,y)  =\int_X f(x,b)d\nu(b)  \hskip 1cm T^{-1}_{\mathcal G}( f)(x,y)  =\int_X f(b,y)d\nu(b)   $$
$$   \widehat T'_{\mathcal G}(\rho(f))(x,y)  = f(y,y)   $$

\vskip 1cm
Let's recall what is a $\mathcal G$-space or an action of $\mathcal G$ (\cite{AR} chap.2), 
\subsubsection{\bf{Notation}}
\label{pringuet}
Let $X,Y$ be two  Borel  spaces, let us call   {\bf fibration} of $X$ by $Y$ any Borel map $\flat: X \to Y$ which is onto. When $Y = \mathcal G^0$ and for any $i \in \{s, r  \}$, let     $X_\flat \times _i \mathcal G $ be the fiber product of $X$ and $\mathcal G$ which is the set   $\{ (x,g) \in X \times \mathcal G / \flat(x) = i(g) \} $ .

\subsubsection{\bf{Definition}} 
\label{iep}
A (right)  $\mathcal G$-space is a Borel space $X$ endowed with a fibration   $\flat: X \to \mathcal G^0$ and a Borel  map $(x,g) \mapsto x.g $  from $ X _\flat \times_r   \mathcal G$ to $X$ such that:

i) For all $(x,g) \in X _\flat \times_r   \mathcal G$, $\flat(x.g) = s(g)$ and $ x.\flat(x) = x$

ii) For all $(x,g_1)  \in X _\flat \times_r   \mathcal G$ and all $g_2 \in \mathcal G^{s(g_1)}$   then  $x.(g_1g_2)  = (x.g_1).g_2$

\subsubsection{\bf{Definition}} 
A (right)  locally compact $\mathcal G$-space is a locally compact space $X$ endowed with a structure of   (right)  $\mathcal G$-space such that $\flat $ is open and continuous and $(x,g) \mapsto x.g $ is continuous.

\vskip 0.1cm
\vskip 1cm

One can also say that $\mathcal G$ acts on $X$. Let us now suppose that $X$ is endowed with a (positive faithful) Radon measure $\theta$ such that $\flat_\star\theta$  is absolutely continuous w.r.t. $\nu$,  then one easily sees that $b: L^\infty(\mathcal G^0,\nu) \to L^\infty(X,\theta)$ defined by $b(f) = f \circ \flat$ is a *(anti)isomorphism,   also for any $i$ in $\{s,r\}$  the von Neumann algebra  $ L^\infty(X,\theta) \surl{\nonumber_{b} \star_{i_{\mathcal G}}}_{L^\infty(\mathcal G^0,\nu)}  L^\infty(\mathcal G ,\mu) $ is canonically isomorphic to  $    L^\infty(X _\flat \times_i  \mathcal G ,\theta _\flat \times_i \mu) $ for a suitable Radon measure $\theta _\flat \times_i \mu$ on $X _\flat \times_i   \mathcal G$,  one easily proves that:

\subsubsection{\bf{Proposition}}
\label{action}
{\it  Let  $X$ be a (right)   locally compact  $\mathcal G$-space endowed with a Radon measure $\theta$ such that $\flat_\star\theta$  is absolutely continuous w.r.t. $\nu$,  let $b: L^\infty(\mathcal G^0,\mu) \to  L^\infty(X,\theta) $ be defined by  $b(f) = f \circ \flat$ for any $f$ in $L^\infty(X,\theta) $, let  $\mathfrak a: L^\infty(X,\theta) \to   L^\infty(X _\flat \times_r   \mathcal G ,\theta _\flat \times_r \mu) $ be  defined by $\mathfrak a(f)(x,g) = f(x.g)$, then $(b,\mathfrak a)$ is an action of $\mathfrak G(\mathcal G)$ on $L^\infty(X,\theta)$.}

\vskip 0.9cm

In the conditions above, let's  suppose $X$ is second countable, then  the crossed product $L^\infty(X,\theta)  \rtimes_{\mathfrak a} \mathfrak G(\mathcal G)$ is also   $\widehat{\mathfrak G}'( X_\flat \times  _{r  }\mathcal G)$,  where $X_\flat \times  _{r  }\mathcal G $  is given its structure of semi direct product.    The existence of $b$, given by Proposition \ref{action},  leads to   integral decompositions $\theta= \int_{\mathcal G^0}   \theta^u d\nu(u)$,   $L^2(X,\theta) =  \int_{\mathcal G^0}   L^2(X^u, \theta^u) d\nu(u)$ and $L^\infty(X,\theta) = \int_{\mathcal G^0}   L^\infty(X^u) d\nu(u)$, where $X^u = \{x \in X/ \flat(x) = u\}$ and the support of $\theta^u$ is $X^u$ for any  $u \in \mathcal G^0$. 

For any Radon measure $\theta'$ on $X$, we shall say that $\theta'$ is invariant under  $\mathfrak a$ if and only if, for any $g\in \mathcal G$, and any $f \in \mathcal K( X^{s(g)})     $ one has: $\int f(x) d\theta^{s(g)} =  \int f(y.g)d\theta^{r(g)} $. 

\subsubsection{\bf{Lemma}}
\label{vergn}
{\it If $\theta$ and $\theta'$ are two  Radon  measures on $X$,  invariant  under  $\mathfrak a$, then  the Radon Nikodym  derivative $\frac{d\theta'}{d\theta}$ is an element of $L^\infty(X,\theta)^{\mathfrak a}$.}
\vskip0.1cm
\begin{dm} This is easy and a basic consequence of    \cite{E4} 7.5 to 7.8.
\end{dm}

\vskip 1cm

Let us now give a description of the crossed product $L^\infty(X,\theta)  \rtimes_\alpha \mathfrak G(\mathcal G)$   using a certain  $*$-algebra representation.
The vector space $\mathcal K (X_\flat \times _{r } \mathcal G )$ can be given a    $*$-algebra structure denoted by $(\mathcal K (X_\flat \times  _{r }\mathcal G ),  \star,\# )$. 

For any $F,F'$ in  $\mathcal K(X_\flat \times  _{r  }\mathcal G  )$ and any  $(x,g)$ in $X_\flat \times _r\mathcal G$, one has:
 $$ F   \star F'(x,g)  =
 \int  F(x,h)F'(x.h,  h^{-1}g)d\lambda^{r(g)}(h)  $$
 
 $$ F^{\#}(x,g)  =
\overline{F(x.g, g^{-1}})\delta(g^{-1})  \  \  $$

One can define a    representation of $\mathcal K(X_b \times _r\mathcal G )$ in $L^2(X_b \times _r\mathcal G, \theta _\flat \times_r \mu)$, let us note it    $\mathfrak R$, it is defined for any $\xi$ in  $L^2(X_b \times _r\mathcal G, \theta _\flat \times_r \mu )$,  any $F$ in $\mathcal K(X_b \times _r\mathcal G )$ and  $\theta _\flat \times_r \mu$-almost any $(x,g)$ in  $X_b \times _r\mathcal G$ by: 
 
  $$    \mathfrak R (F)\xi(x,g) =  \int F(x.g, h) \delta(h)^{-\frac{1}{2}} \xi(x, gh)  d\lambda^{s(g)}(h) $$

The crossed product $L^\infty(X,\theta)  \rtimes_{\mathfrak a } \mathfrak G(\mathcal G)$ is generated by the image of $\mathfrak R$. More precisely, for any $f \in \mathcal K(X )$ and any $k \in \mathcal K(\mathcal G)$, one has: $\mathfrak a(f)(1 _b\otimes_r \rho(k)) = \mathfrak R(f \otimes k)$.

\vskip 1.5cm

\section{A generalization  of the matched pair procedure}

\vskip 1cm
\subsection{Measured matched   pair of groupoids }
\label{matched}
Now let's explain a triple    extension of   the commutative examples (in \ref{R}), those  studied in \cite{Val4} in finite dimension, and in \cite{VV} (or   \cite{BSV}) in the   quantum groups case.  

All the groupoids involved  will be  as mentioned in \ref{R}  till the end; mimicking the case of matched pairs of locally compact groups, let's give the following definition:

\subsubsection{\bf{Definition}}
\label{mapa}
{\it  Let $\mathcal G$, $\mathcal G_1$, $\mathcal G_2$  be  measured groupoid such that $ \mathcal G_1, \mathcal G_2$  are      closed subgroupoids of  $\mathcal G$.  We shall say that $\mathcal G_1, \mathcal G_2$ is a matched pair of measured subgroupoids of $\mathcal G$ if  and only if:

i) $\mathcal G_1 \cap  \mathcal G_2 = \mathcal G^0$

ii)     $ \mathcal G_1  \mathcal G_2 :=   \{g_1g_2/ g_1 \in \mathcal G_1, g_2 \in {\mathcal G}_{2}^{s(g_1)}\}$    is $\mu$-conegligible in $\mathcal G$

iii) the  measure $\nu$ on $\mathcal G^0$  is quasi invariant  for the three Haar systems.  }

\vskip 0.5cm

\subsubsection{\bf Remark}
\label{doubout971}
The condition iii) is absolutely necessary to obtain, as in the case of groups (\cite{BSV} prop. 3.2), the Haar system of $\mathcal G$ from those of $\mathcal G_1$ and $\mathcal G_2$.  for instance if $\mathcal G$ is a principal groupoid of the form $X\times X$ where $X$ is any locally compact space, let's choose $\mathcal G_1 = \mathcal G $ and $\mathcal G_2 = \mathcal G^0 = X$, then if $\nu$ and $\nu_1$ and $\nu_3$ are any Radon measures on $X$, one can construct the Haar systems $(\delta_x \times \nu)_{x \in X}$ and $(\delta_x \otimes \nu_1)_{x \in X}$ on $\mathcal G$ and $\mathcal G_1$ and the quasi invariant measures $\nu $ and $\nu_1$ respectively, there is no hope to give any formula connecting the Haar systems.  We shall now prove such a formula, when the   groupoids are given the same quasi invariant measure,  using an argument similar to that of  \cite{BSV} prop. 3.2 and the  uniqueness  condition of M.Enock (\cite{E4} Corollary 7.8) recalled in \ref{vergn}.

 \vskip 0.5cm
\subsubsection{\bf Lemma}
\label{martyr}
{\it 1) The von Neumann fiber product $L^\infty(\mathcal G_1,\mu_1)\surl{\nonumber_{s_1} \star_{s_2}}_{L^\infty(\mathcal G^0,\nu)}L^\infty(\mathcal G_2,\mu_2)$ is isomorphic to $  L^\infty(\mathcal G_1{_s \times_s} \mathcal G_2,\mu_1{_s \times_s}\mu_2)$, where the measure $\mu_1{_s \times_s}\mu_2$ is given for any $f \in \mathcal K( G_1{_s \times_s} \mathcal G_2)$ by:
$$(\mu_1{_s \times_s}\mu_2)(f) =  \int_{\mathcal G^0}\int_{G_1{_r \times_r} \mathcal G_2}f(g_1^{-1},g_2^{-1})\delta_{1}(g_1^{-1})\delta_2(g_2^{-1})d\lambda_1^u(g_1)d\lambda_2^u(g_2)d\nu(u)$$

2) The restriction to $  L^\infty(\mathcal G_1{_s \times_s} \mathcal G_2,\mu_1{_s \times_s}\mu_2)$ of the natural action by left multiplication of $\mathfrak G(\mathcal G_1{ \times } \mathcal G_2)$ on $L^\infty(\mathcal G_1 \times \mathcal G_2, \mu_1 \times \mu_2)$ is well defined;  if  $(\mathfrak a, r_1{_s \otimes_s} r_2)$ is  this action, then for any  $f\in L^\infty(\mathcal G_1{_s \times_s} \mathcal G_2,\mu_1{_s \times_s}\mu_2)$,  $\mu_1{_s \times_s}\mu_2$-any element $(g_1,g_2) \in G_1{_s \times_s} \mathcal G_2$, $\mu_1$-any $h_1 \in \mathcal G_1$, $\mu_2$-any $h_2 \in \mathcal G_2$ ,one has:
$$\mathfrak a(f)(g_1,g_2,h_1,h_2) = f(h_1^{-1}g_1,h_2^{-1}g_2)$$
and  $\mu_1{_s \times_s}\mu_2$ is invariant under $\mathfrak a$ in the sense of the end of paragraph \ref{R}.}

\vskip 0.1cm

\begin{dm}
The assertion 1)  is a simple calculation (see  for instance \cite{Val4} 3.1) and the second one,   an obvious consequence of 1), as one deals with the left multiplication of $\mathcal G_1 \times \mathcal G_2$ and its canonical Haar system. 
\end{dm}

\vskip 0.5cm

\subsubsection{\bf Lemma}
\label{972}
{\it  Let $\tilde \mu$ be the measure on  $ \mathcal G_1{_s \times_s} \mathcal G_2$ defined for any $f \in \mathcal K( \mathcal G_1{_s \times_s} \mathcal G_2)$ by:
$$ \tilde \mu(f) = \int_{u \in \mathcal G^0} \int_{\mathcal G_1\mathcal G_2} \tilde f(\theta^{-1}(g))    d\lambda ^u(g )d\nu(u) $$
where $\theta: \mathcal G_1{_s \times_s} \mathcal G_2 \to \mathcal G_1 \mathcal G_2$ is given by $\theta(g_1,g_2)= g_1g_2^{-1}$ and $\tilde f(g_1,g_2) = f(g_1,g_2)\delta(g_2)$. 
Then $\tilde \mu$ is invariant under  the action  $\mathfrak a $.
}
\vskip 0.1cm

\begin{dm} This is the same argument as in lemma 4.10 of \cite{VV}.
 \end{dm}

\subsubsection{\bf{Proposition}}
\label{use}
{\it Let $\mathcal G_1, \mathcal G_2$ be  a matched pair of measured subgroupoids of $\mathcal G$, and let $\delta_i$  be   the modular function of $\nu$ (relatively to $\mathcal G_i$) for $i= 1,2$, then 
 up to normalization of the Haar systems,  for any  $u \in \mathcal G^0$ and any $f \in \mathcal K(\mathcal G)$, one has: 
 \begin{align*}
\int fd\lambda^u &=   \iint_{\mathcal G_1\times \mathcal G_2} f({g}_1{g}_2)\delta({g}_2))\delta_2({g}^{-1}_2) d\lambda_2^{s(g_1)}(g_2)d\lambda^u_1 (g_1)
\\
&=   \iint_{\mathcal G_2\times \mathcal G_1} f({g}_2g_1)\delta(g_1)\delta_1({g}^{-1}_1) d\lambda_1^{s(g_2)}(g_1)d\lambda^u_2(g_2) 
\end{align*}  }
\vskip 0.2cm
\begin{dm}
Due to the lemmas \ref{martyr} and \ref{972}, one can apply  lemma   \ref{vergn} to $\mu_1{_s \times_s}\mu_2$ and $\tilde\mu$, so there exists a function $h \in L^\infty(\mathcal G_1{_s \times_s} \mathcal G_2,\mu_1{_s \times_s}\mu_2)^\mathfrak a$ such that $\tilde \mu = h(\mu_1{_s \times_s}\mu_2)$, hence   for  $\mu_1{_s \times_s}\mu_2$-any  $(g_1,g_2) \in \mathcal G_1{_s \times_s} \mathcal G_2$: $h(g_1,g_2) = h(s(g_1),s(g_2))$. For any  $F \in \mathcal K(\mathcal G)$, if $f \in    L^\infty(\mathcal G_1{_s \times_s} \mathcal G_2,\mu_1{_s \times_s}\mu_2)$  is defined   for any $(g_1,g_2) \in \mathcal G_1{_s \times_s} \mathcal G_2$ by $f(g_1,g_2) = F(g_1g_2^{-1})\delta(g_2^{-1})$ then $\tilde f \circ \theta = F$, so one has:

\begin{align*}
 \int_{u \in \mathcal G^0} &\int_{\mathcal G_1\mathcal G_2} F(g) )  d\lambda ^u(g )d\nu(u)
 =  \int_{u \in \mathcal G^0} \int_{\mathcal G_1\mathcal G_2} \tilde f(\theta^{-1}(g))    d\lambda ^u(g )d\nu(u) 
=  \tilde \mu(f) \\
&=  (\mu_1{_s \times_s}\mu_2)(hf) \\
 &= \int_{\mathcal G^0}\int_{G_1{_r \times_r} \mathcal G_2}h(g_1^{-1},g_2^{-1})f(g_1^{-1},g_2^{-1})\delta_{1}(g_1^{-1})\delta_2(g_2^{-1})d\lambda_1^u(g_1)d\lambda_2^u(g_2)d\nu(u)\\
 &= \int_{\mathcal G^0} \int_{G_1{_r \times_r} \mathcal G_2} h(r(g_1),r(g_1))F(g_1^{-1}g_2)\delta_{1}(g_1^{-1})\delta(g_2)\delta_2(g_2^{-1})d\lambda_2^u(g_2)d\lambda_1^u(g_1)d\nu(u)\\
 &= \int_{\mathcal G_1}\big( \int_{G_2} h(r(g_1),r(g_1))F(g_1^{-1}g_2)\delta_{1}(g_1^{-1})\delta(g_2)\delta_2(g_2^{-1}) d\lambda_2^{r(g_1)}(g_2) \big ) d\mu_1(g_1)\\
 &= \int_{\mathcal G_1}\big( \int_{G_2} h(s(g_1),s(g_1))F(g_1 g_2) \delta(g_2)\delta_2(g_2^{-1}) d\lambda_2^{s(g_1)}(g_2) \big ) d\mu_1(g_1)\\
 &= \int_{\mathcal G^0} \int_{G_1{_r \times_r} \mathcal G_2} F(g_1g_2) \delta(g_2)\delta_2(g_2^{-1})\lambda_2^{s(g_1)}(g_2) h(s(g_1),s(g_1)) d\lambda^u_1(g_1) d\nu(u)\end{align*}
 
This gives the first equality of the Proposition, if one replaces  the Haar system $\lambda_1^u$ by $k\lambda_1^u$, where $k$ is  defined by $k(g_1) = h(s(g_1),s(g_1))$ which is still a Haar system (as $k(g_1)$ depends only on $s(g_1)$). The second equality is proven a similar way.
\end{dm}

\subsubsection{\bf{Remark}}
\label{vengeance1}

 If $\mathcal G_1,\mathcal G_2$ is a measured matched pair, as $\mathcal G_2\mathcal G_1 = (\mathcal G_1\mathcal G_2)^{-1}$ and $\mathcal G_1\mathcal G_2 \cap \mathcal G_2\mathcal G_1$ is conegligible then  $\mathcal G_2,\mathcal G_1$ is also a measured matched pair.

 \vskip  1cm 
\subsection{Families of examples}

\subsubsection{ \bf Example}
If $\mathcal G$ is a group then a matched pair of groupoids is a matched pair of groups as  in   \cite{BSV}.

\subsubsection{ \bf Example}
If $\mathcal G$ is finite then a  matched pair is exactly a matched pair as defined in    \cite{Val2}.

\subsubsection{ \bf Example}
If $\mathcal G$ is a group bundle (i.e. $s=r$),  as for any $u \in \mathcal G^0$,  $\mathcal G^u = \mathcal G_u^u$ is a group, the  matched pairs are exactly  group bundles over $G^0$ such that for any $u \in \mathcal G^0$,  $\mathcal G_1^u, \mathcal G^u_2$ is a matched pair of groups in $\mathcal G^u$. The Haar systems being  continuous families of Haar measures, in that case, proposition \ref{use} is a completely direct consequence of  \cite{BSV} prop. 3.2.
 
\subsubsection{\bf Example}
\label{vengeance2}
If $\mathcal G$ is a locally compact transformation groupoid of the form $X \times G$,  where $G$ is a  group matched pair  $G_1G_2$  in the sense of \cite{BSV}, acting on the right on $X$.  Let $\mathcal G_1$ be equal to $X \times G_1$, $\mathcal G_2$ be  equal to  $X \times  G_2 $.  One can  consider  the canonical   Haar systems $(\delta_x \times ds)_{x \in X}$, $(\delta_x \times ds_1)_{x \in X}$ $(\delta_x \times ds_2)_{x \in X}$ associated with the Haar measures of the groups, hence $\mathcal G_1$, $\mathcal G_2$ is clearly a matched pair.

 Let $\nu$ be  a quasi invariant measure on $\mathcal G^0 = X$ w.r.t.  the action of $G$,  let $\rho$  be the Radon Nikodym   cocycle for $\nu$ and the action, this means that  for any $g \in G$ and $h \in \mathcal K(X)$ one has $\int h(x.g) d\nu(x) = \int \rho(x,g)h(x)d\nu(x)$. Let $\Delta$ (resp. $\Delta_1$, resp.  $\Delta_2$) be the modular function of the group $G$ (resp.  $G_1$, resp. $G_2$), hence by \cite{R} 3.21 one has:
$\delta(x,g) = \frac{\Delta(g)}{\rho(x,g)}$ (resp $\delta_i(x,g) = \frac{\Delta_i(g)}{\rho(x,g)}$ for $i=1,2$). Therefore  due to \cite{BSV} prop. 3.2, for any $f \in \mathcal K (X \times \mathcal G)$ and any $x \in X$, one has:

\begin{align*}
\delta_x& \times ds(f)
= \int f(x,g)ds(g) = \int_{G_1}\int_{G_2}f(x,g_1g_2)\Delta(g_2)\Delta_2(g_2)^{-1}ds_2(g_2)ds_1(g_1) \\
&= \int_{G_1}\int_{G_2}f((x,g_1)(xg_1,g_2)) \frac{\Delta(g_2)}{\rho(xg_1 ,g_2)}(\frac{\Delta_2(g_2)} {\rho(xg_1,g_2)})^{-1}ds_2(g_2)ds_1(g_1)  \\
&= \int_{G_1}\int_{G_2}f((x,g_1)(xg_1,g_2)) \delta( xg_1,g_2) \delta_2(xg_1,g_2)^{-1}ds_2(g_2)ds_1(g_1)   \\
&= \int_{\mathcal G_1}\int_{\mathcal G_2}f((y,g_1)(z,g_2)) \delta( z,g_2) \delta_2(z,g_2)^{-1}(\delta_{s(x,g_1)}\times ds_2)(z,g_2)(\delta_x \times ds_1)(y,g_1) 
\end{align*}
which gives Proposition \ref{use} in that case.

\subsection{The case of  principal and transitive groupoids }
\label{vengeance3}
 Let's describe matched pairs in the case where $\mathcal G$ is principal proper and transitive.
 \vskip0.5cm
  So    we shall suppose that $\mathcal G = X \times X$  where $X$ is an Hausdorff  locally compact space together with some Radon  measure $\nu$, and $\mathcal G$ is given its natural locally compact groupoid structure, Haar system and quasi invariant measure $\nu$, as it was explained in remark  \ref{doubout971}. 
 
 Let us describe a family of examples and let us show  that all  matched pairs of a  principal proper and transitive groupoid  are of this type  when X is compact.
 
 \vskip 0.2cm
 
 We shall suppose that $X$ is equal to $X_1 \times X_2$ where $X_1$ and $X_2$ are two Hausdorff locally compact spaces and $\nu = \nu_1 \times \nu_2$, where $\nu_i$ is a Radon measure on $X_i$ for $i = 1,2$. Let $\mathcal R_1$ and $\mathcal R_2$ be the  equivalence relations associated with the natural projections,  so one has:
 $$ \forall (a,b), (c,d) \in X:  \  \  (a,b) \mathcal R_1 (c,d) \  \  \mathrm{iff}   \  \  a= c$$
 $$ \forall (a,b), (c,d) \in X:  \  \   (a,b) \mathcal R_2 (c,d)  \  \  \mathrm{iff}   \  \   b= d$$
 Let $\mathcal G_i$ be the sub groupoids of $\mathcal G = X \times X$ associated with $\mathcal R_i$ for $i= 1,2$. One has:
 $$ \mathcal G_1 = \underset{x_1 \in X_1}\sqcup   \  \ \{x_1\}\times X_2\times  \{x_1\}\times X_2$$
  $$ \mathcal G_2 = \underset{x_2 \in X_2}\sqcup   \  \ \    X_1 \times \{x_2\} \times X_1 \times\{x_2\}$$
 
 As $X$ is Hausdorff, clearly $\mathcal G_1$ and $\mathcal G_2$ are closed in $\mathcal G$. For any $(x_1,x_2) \in X (=\mathcal G^0)$, let $\lambda^{(x1,x2)}$ be equal to $\delta_{(x_1,x_2)}\times \nu$ which gives the canonical Haar system of $\mathcal G$ and let's define two other Radon  measures  on $\mathcal G_1$ and $\mathcal G_2$ respectively by the formulas:
 $$\lambda^{(x1,x2)}_1 = \delta_{(x_1,x_2)} \times \delta_{x_1}\times {\nu_2}$$
 $$\lambda^{(x1,x2)}_2 = \delta_{(x_1,x_2)} \times {\nu_1}\times \delta_{x_2}$$
  
 \subsubsection{\bf Lemma}
 \label{ced}
 {\it   $(\mathcal G_1, (\lambda^u_1)_{u \in X},\nu)$, $(\mathcal G_2, (\lambda^u_2)_{u \in X},\nu)$ is a   matched pair in $(\mathcal G, (\lambda^u)_{u \in X},\nu)$.}
 \vskip 0.3cm
 \begin{dm}
 For any $(a,b,c,d) \in \mathcal G$, one has: $(a,b,c,d) = (a,b,a,d).(a,d,c,d)$ so $\mathcal G = \mathcal G_1 \mathcal G_2$,  and $\mathcal G_1 \cap \mathcal G_2 = \underset{(x_1,x_2) \in X}\sqcup   \  \ \{x_1\}\times\{x_2\}\times  \{x_1\}\times \{x_2\}= \mathcal G^0$  .
   
 For any $(x_1,x_2) \in X$, the support of $ \lambda^{(x_1,x_2)}_1$ is clearly $\mathcal G^{(x_1,x_2)}_1$. For any $(x_1,x_2,x_1,z_2)$ in $\mathcal G_1$ and for any $f \in \mathcal K(\mathcal G_1)$, one has:
 \begin{align*}
 \int_{\mathcal G_1}f((x_1,x_2,x_1,z_2)t) d\lambda^{(x_1,z_2)}_1(t)
 &=  \int_{X_1 \times X_2}f(x_1,x_2,t_1,t_2) \delta_{x_1}(t_1)d\nu_2(t_2) \\
 &=  \int_{  X_2}f(x_1,x_2,x_1,t_2) d\nu_2(t_2) \\
 &=  \int_{  \mathcal G_1}f(t) d\lambda^{(x_1,z_2)}_1(t)
  \end{align*}
 
 One easily deduces that $(\lambda^u_1)_{u \in X}$ is a continuous Haar system for $\mathcal G_1$.
 
  Let $\mu_1$ be equal to $\int_X  \lambda^u_1d\nu(u)$,  for any  $f$  in $\mathcal K(\mathcal G_1)$   one has:
  \begin{align*}
  \int_{\mathcal G_1} f(z^{-1})d\mu_1(z)
  &= \int_X\int_{\mathcal G_1}f((y_1,y_2,z_1,z_2)^{-1})\delta_{(x_1,x_2)} \times \delta_{x_1}\times {d\nu_2}(y_1,y_2,z_1,z_2)d\nu(x_1,x_2)\\
  &=  \int_X\int_{X_2}f((x_1,x_2,x_1,z_2)^{-1}) d\nu_2 ( z_2)d\nu(x_1,x_2)\\
  &=  \int_{X_1}\int_{X_2}\int_{X_2}f(x_1,z_2,x_1,x_2) d\nu_2 ( z_2)d\nu_1(x_1)d\nu_2(x_2)
  \end{align*}
  Hence using Fubini's theorem, one has : $ \int_{\mathcal G_1} f(z^{-1})d\mu_1(z) = \int_{\mathcal G_1} f(z)d\mu_1(z)$, so $\nu$ is quasi invariant, even invariant, relatively to the Haar system $(\lambda^u_1)_{u \in X}$ and  $(\mathcal G_1, (\lambda^u_1)_{u \in X},\nu)$ is a measured groupoid. For similar reasons, $(\mathcal G_2, (\lambda^u_2)_{u \in X},\nu)$ is also a measured groupoid and $\nu$ is invariant for the three measured groupoids, so $(\mathcal G_1, (\lambda^u_1)_{u \in X},\nu)$, $(\mathcal G_2, (\lambda^u_2)_{u \in X},\nu)$ is a matched pair.

 \end{dm}

  \vskip 0.5cm
Let us verify proposition \ref{use} in that case, for any $h \in \mathcal K(\mathcal G)$ and any $(x_1,x_2)\in X$, one has:
  
  \begin{align*}
 & \int_{\mathcal G_2}  \int_{\mathcal G_1} h((a_1,a_2,b_1,b_2)(c_1,c_2,d_1,d_2))d\lambda^{s(a_1,a_2,b_1,b_2)}_2(c_1,c_2,d_1,d_2)d\lambda^{(x_1,x_2)}_1(a_1,a_2,b_1,b_2)\\
 &=  \int_{ X_2}  \int_{\mathcal G_1} h((x_1,x_2,x_1,b_2)(c_1,c_2,d_1,d_2))d\lambda^{(x_1,b_2)}_2(c_1,c_2,d_1,d_2)d\nu_2(b_2)\\
 &=  \int_{ X_2}  \int_{ X_1} h(x_1,x_2,d_1,b_2)d\nu_1(d_1)d\nu_2(b_2)
  =  \int_{ X_1}  \int_{ X_2} h(x_1,x_2,d_1,b_2)d\nu_1(d_1)d\nu_2(b_2)\\
  &= \int_{\mathcal G} h(z)d\lambda^{(x_1,x_2)}(z)
  \end{align*}
  
  This proves that proposition \ref{use} is here a reformulation of Fubini's theorem.

 \vskip 1cm

Finally, let's prove  that we have described all possible examples when $X$ is Hausdorff and compact. So we suppose given a matched pair $\mathcal G_1$, $\mathcal G_2$ in $X\times X$,  there exist two  equivalence relations $\mathcal R_1$ and $\mathcal R_2$ associated to two partitions $(X^1_\alpha), (X^2_\beta)$ of $X$, using classical arguments (see \cite{Go} chap 1 par.4), $\mathcal R_1$ and $\mathcal R_2$ are closed and Hausdorff, so   each element of the partitions is  closed as  a subset of $X$,  hence  compact in $X$. As well-known, due to the fact that $\nu$ is quasi invariant for $\mathcal G_i$ ($i= 1,2$), there exist  Radon measures $\Lambda_i$ on $X / \mathcal R_i$ and Borel functions $h_i : X \to \mathbb R^+_\star$ such that if one denotes $\alpha_i : X \to X / \mathcal R_i$  the usual projection, one has $\mu = h_i(\Lambda_i \circ \alpha_i)$.

 \subsubsection{\bf Proposition}
 \label{allegre}
 {\it The map   \  $\alpha_1 \times \alpha_2: \ \ X \to X / \mathcal R_1 \times X / \mathcal R_2$ defined for any $x \in X$ by:   $(\alpha_1 \times \alpha_2)(x) = (\alpha_1(x), \alpha_2(x))$, realizes a   homeomorphism of compact spaces, and up to normalization $(\alpha_1 \times \alpha_2)(\mu) = \Lambda_1 \otimes \Lambda_2$.}
 \vskip 0.2cm
 \begin{dm}
 The application $\alpha_1 \times \alpha_2$ is clearly continuous, let $x,y \in X$ be such that $(\alpha_1 \times \alpha_2)(x) = (\alpha_1 \times \alpha_2)(y)$ then  $(x,y) \in \mathcal G_1 \cap \mathcal G_2$ so $x=y$, hence $\alpha_1 \times \alpha_2$ is injective. As $\mathcal G_1 \mathcal G_2$ is $\mu$ conegligible in $\mathcal G$  and compact  its complementary is a negligible open set so it is empty and  $\mathcal G_1 \mathcal G_2 =\mathcal G$, this implies   that for  any $x , y  \in X $, there exists  $t \in X$ such that $(x,t)\in \mathcal G_1$ and $(t,y)  \in \mathcal G_2$, hence $\alpha_1(x) =\alpha_1(t)$ and $\alpha_2(y) =\alpha_2(t)$ which means that $(\alpha_1 \times \alpha_2)(t) =  (\alpha_1(x), \alpha_2(y))$, hence  so $\alpha_1 \times \alpha_2$  is onto,   the lemma follows easily. 

 \end{dm}

\vskip 0.6cm

\subsection{ The   mutual actions  of a matched pair of groupoids}
Let's give a generalization to matched pairs of groupoids of the well-known fact that matched pairs of groups act  one on the other.

\vskip 1cm
\subsubsection{\bf Remark}
\label{rollex}
As in \cite{BSV1} Chap 2, if $ (\mathcal G_1,(\lambda_1^u),\nu ), (\mathcal G_2, (\lambda_2^u),\nu )$ is a  given   measured matched pair, then for $i=1$ and  $2$, there exist    Borel functions $p_i: \mathcal G  \to \mathcal G_i$   such that,   for any $g \in \mathcal G_1\mathcal G_2$ , $g = p_1(g)p_2(g)$. Of course, there also exists  two borel almost everywhere defined maps  $p'_i: \mathcal G \to \mathcal G_i$   such that $g = p'_2(g)p'_1(g)$, for any  $g \in \mathcal G_2\mathcal G_1$; so on the $\mu$- conegligible set $ \mathcal G_1\mathcal G_2 \cap  \mathcal G_2\mathcal G_1$, one has: $ g = p_1(g)p_2(g)= p'_2(g)p'_1(g)$.  In this framework, new representations  appear, the {\bf{middle}} ones: 
\vskip 0.6cm

\subsubsection{\bf Lemma and notations}
\label{middle1}
{\it 
  For $\mu$-almost any $g$ in $\mathcal G$,  one has:   $s \circ p_1(g)  = r \circ p_2(g) ,   r \circ p'_1(g) = s\circ p'_2(g)$,  so   there exist two $\mu$-almost everywhere defined maps  such that:  
   $$m=   s \circ p_1  = r \circ p_2  \ \  ,  \  \  \hat m = s \circ p'_2 = r \circ p'_1$$
   let us   note      $m_{\mathcal G}:  f \mapsto f \circ m$ (resp.$\hat m_{\mathcal G}: f \mapsto f \circ m_{\mathcal G}$) the associated representation of  $L^\infty(\mathcal G^0,\nu)$. }
   \vskip 0.1cm
   \begin{dm} As for any $g \in \mathcal G_1\mathcal G_2$, $p_1(g)$ and $p_2(g)$ are composable, the lemma  is obvious. 
   \end{dm}
   \vskip 1cm
\subsubsection{\bf Lemma}
\label{sarcofobe}
{\it  For $i= 1$ or $2$,  let us note  $s_i = s_{\mathcal G_i}$ and $r_i = r_{\mathcal G_i}$. There exists an isomorphism    $U: L^2(\mathcal G,\mu) \to  L^2(\mathcal G_2,\mu_2)\underset{\nu}{_{s_2}\otimes_{r_1}}L^2(\mathcal G_1,\mu_1) (=  L^2(\mathcal G_2  {_{s_2}\times_{r_1}}\mathcal G_1,\mu_2 {_{s_2}\times_{r_1}}\mu_1))$   such that for any $\xi \in \mathcal K(\mathcal G)$ and $\mu_2 {_{s_2}\times_{r_1}}\mu_1$-almost any $(g_2,g_1) \in \mathcal G_2   {_{s_2}\times_{r_1}}\mathcal G_1$, one has:  
$$U\xi(g_2,g_1) = (\frac{\delta}{\delta_{1}})^{\frac{1}{2}}(g_1) \xi(g_2g_1)$$
it  induces   an isomorphism  between   $L^\infty(\mathcal G,\mu)$ and   $L^\infty(\mathcal G_2,\mu_2)\underset{\nu}{_{s_2}\star_{r_1}}L^\infty(\mathcal G_1,\mu_1)$, these two von Neumann algebras   are  also isomorphic to  $ L^\infty(\mathcal G_2    {_{s_2}\times_{r_1}}\mathcal G_1,\mu_2 {_{s_2}\times_{r_1}}\mu_1)$. }
 \vskip 0.2cm
  \begin{dm} This is an obvious consequence of proposition \ref{use} .\end{dm}
  \vskip 1cm

\subsubsection{\bf{Proposition and definition}}
\label{bienfait}
{\it  Let  $\mathfrak a: L^\infty(\mathcal G_2,\mu_2) \to    L^\infty(\mathcal G_2,\mu_2)\underset{\nu}{_{s_2}\star_{r_1}}L^\infty(\mathcal G_1,\mu_1)$ (resp.$\hat{\mathfrak a}: L^\infty(\mathcal G_1,\mu_1) \to  L^\infty(\mathcal G_1,\mu_1)\underset{\nu}{_{r_1}\star_{s_2}}L^\infty(\mathcal G_2,\mu_2)$) be  the   map defined for any $f \in L^\infty(\mathcal G_2,\mu_2)$ (resp.any $h \in L^\infty(\mathcal G_1,\mu_1)$)  and  almost any $(g_2,g_1) \in \mathcal G_2 {_{s_2}\times_{r_1}} \mathcal G_1$  (resp. $(g_1,g_2) \in \mathcal G_1 {_{r_1}\times_{s_2}}  \mathcal G_2$) by:
   $$\mathfrak a(f)(g_2,g_1) = f(p_2(g_2g_1))\  \  \  \mathrm{(resp. } \  \ \hat {\mathfrak  a}(h)(g_1,g_2) =h(p_1(g_2g_1))) $$
   Then the   pair $(s_2,\mathfrak a)$ (resp.$(r_1,\hat{\mathfrak  a})$) is a right (resp.left) action of $\mathfrak G(\mathcal G_1)$ (resp.$\mathfrak G(\mathcal G_2)$) on $L^\infty(\mathcal G_2,\mu_2)$ (resp.$L^\infty(\mathcal G_1,\mu_1)$), moreover one has: $\mathfrak a \circ r_2   = \hat{\mathfrak a} \circ s_1= m$ and $\mathfrak a \circ s_2 = s$.    }
\vskip 0.7cm
\begin{dm}  For all $\phi  \in L^\infty(\mathcal G^0,\nu)$, and $\mu_2 {_{s_2}\times_{r_1}}\mu_1$-almost any $(g_2,g_1) \in \mathcal G_2\times \mathcal G_1$ one has: 
\begin{align*}
\mathfrak a(s_2(\phi))(g_2,g_1) =  s_2(\phi))(p_2(g_2g_1) )=  \phi(s( p_2(g_2g_1)) ) = \phi(s(  g_1) ) =  (1 {_{s_2}\otimes}_{r_1}s_1(\phi))(g_2,g_1 )
\end{align*}
so $\mathfrak a(s_2(\phi))  =  1 {_{s_2}\otimes}_{r_1}s_1(\phi)$, up to the identification of \ref{sarcofobe}, the relation  $\mathfrak a \circ r_2   = m$ is obtained in the same way.

Let $f$ be any element of $  L^\infty(\mathcal G_2,\mu_2)$ and     let  $(g_2,g_1,h_1)$ be any element  of $  \mathcal G_2 \times  \mathcal G_1 \times  \mathcal G_1$ such that   such that $g_2g_1$ and  $g_2g_1h_1$ exist and  are in $ \mathcal G_1\mathcal G_2$.  On the one hand, we have:
\begin{align*}
(\mathfrak a  \surl{\nonumber_{s_2} \star_{r_1}}_{L^\infty(\mathcal G^0,\nu)} i)\mathfrak a(f)(g_2,g_1,h_1) = \mathfrak a(f)(p_2(g_2g_1),h_1))=  f(p_2(p_2(g_2g_1)h_1)) 
\end{align*}
on the other hand:
\begin{align*}
(i \surl{\nonumber_{s_2} \star_{r_1}}_{L^\infty(\mathcal G^0,\nu)}\Gamma_1) \mathfrak a (f)(g_2,g_1,h_1) = \mathfrak a(f)(g_2, g_1h_1))=  f(p_2(g_2g_1h_1)) 
\end{align*}

But one has  $p_2(p_2(g_2g_1)h_1) = p_2(p_1(g_2g_1)p_2(g_2g_1)h_1) = p_2(g_2g_1h_1)$, hence:
$$ (\mathfrak a  \surl{\nonumber_{s_2} \star_{r_1}}_{L^\infty(\mathcal G^0,\nu)} i)\mathfrak a  =     (i \surl{\nonumber_{s_2} \star_{r_1}}_{L^\infty(\mathcal G^0,\nu)}\Gamma_1) \mathfrak a$$
The demonstration for $\hat{\mathfrak a}$ is quite similar.
\end{dm}
 \subsubsection{\bf{Remark}}
\label{dublin}
The crossed product $L^\infty(\mathcal G_2,\mu_2) \rtimes_{\mathfrak a} \mathfrak G(\mathcal G_1)$ is, up to an isomorphism,   the image in $\mathcal L(L^2(\mathcal G, \mu))$,   of  the map   $\mathcal R$,   defined for any $F \in \mathcal K (\mathcal  G_2   {_{s_2}\times_{r_1}} \mathcal G_1 )$ and $\mu$-almost any $g \in \mathcal G$ by:
$$ \mathcal R(F)\xi(g) = \int_{\mathcal G_1}F(p_2(g),g'_1)\delta(g'_1)^{-\frac{1}{2}} \xi(gg'_1)d\lambda^{s(g)}_1(g'_1)$$
\subsubsection{\bf{Remark}}
\label{modeste} 
Proposition \ref{bienfait} generalizes the fact,  proven in \cite{VV} chap.4, that if $G_1,G_2$ is a matched pair of groups, then  there exists a canonical   action  of $G_1$ on $L^\infty(G_2)$ (resp. $G_2$ on $L^\infty(G_1)$) coming from a map    $\beta$  (resp. $\alpha$)  such that,  up to negligible sets, for any $g_i \in  G_i$ ($i=1,2$), one has $g_1g_2^{-1} = \beta_{g_1}(g_2)^{-1}\alpha_{g_2}(g_1)$.

\vskip 0.9cm

\subsection{ A pseudo multiplicative unitary associated with a matched pair  }
\subsubsection{\bf{Lemma}}
   \label{middle}
 {\it For any $f,f' \in \mathcal K(\mathcal G)$, $\nu$-almost any $u \in \mathcal G^0$, one has:
   $$ \frac{d\omega_{f ,f'}\circ m_{\mathcal G}}{d\nu}(u) = \int_{\mathcal G_1\times\mathcal G_2 }f\overline{f'}(g_1^{-1}g_2)\delta_1(g_1^{-1})\delta(g_2)\delta_2(g_2^{-1})d\lambda_1^u(g_1)d\lambda_2^u(g_2)  $$}
     \newline
  \begin{dm} 
  For any $f,f' \in \mathcal K(\mathcal G)$,  $h \in \mathcal K(\mathcal G^{0})$, and $\nu$-almost any $u \in \mathcal G^0$, one has: 
  $$\begin{array}{lll}
&(\omega_{f ,f'}\circ m_{\mathcal G})(h)
=
 \\

\\
&=\displaystyle{\int_{\mathcal G^0}\int_{\mathcal G}h(m(g))(f\overline{f'})(g)d\lambda^u(g)d\nu(u)}
 
\\

\\
&=\displaystyle{\int_{\mathcal G^0}\int_{\mathcal G_1\times \mathcal G_2} h(m(g_1g_2))(f\overline{f'})({g}_1{g}_2)(\delta\delta_2^{-1})({g }_2) d\lambda_2^{s(g_1)}(g_2)d\lambda_1^u(g_1)d\nu(u)} &  \mathrm{}\, \textit{by \ref{use}}
\\

\\
&=\displaystyle{\int_{\mathcal G^0}\int_{\mathcal G_1\times \mathcal G_2} h(s(g_1))(f\overline{f'})({g}_1{g}_2)(\delta\delta_2^{-1})({g }_2) d\lambda_2^{s(g_1)}(g_2)d\lambda_1^u(g_1)d\nu(u)} & \mathrm{ }\, \textit{ } 
\\

\\
&=\displaystyle{\int_{\mathcal G_1} h(s(g_1))\int_{\mathcal G_2}(f\overline{f'})({g}_1{g}_2)(\delta\delta_2^{-1})({g }_2) d\lambda_2^{s(g_1)}(g_2)d\mu_1(g_1)}

\\

\\
&=\displaystyle{\int_{\mathcal G_1} h(r(g_1))\delta_1(g_1^{-1})\int_{\mathcal G_2}(f\overline{f'})({g}_1^{-1}{g}_2)(\delta\delta_2^{-1})({g }_2)d\lambda_2^{r(g_1)}(g_2)d\mu_1(g_1)} &  \mathrm{ }\, \textit{} 

\\

\\
&  \hskip 4cm   \mathrm{by \ \ making \ \ the \ \ change \ \ of \ \ variable: g_1 \mapsto g_1^{-1}} &  \mathrm{}\, \textit{} 

\\

\\
&=\displaystyle{\int_{\mathcal G^0}h(u)\int_{\mathcal G_1}\int_{\mathcal G_2}(f\overline{f'})({g}_1^{-1}{g}_2)\delta_1(g_1^{-1})(\delta\delta_2^{-1})({g }_2) d\lambda_2^{r(g_1)}(g_2)d\lambda_1^{u}(g_1)d\nu(u)} &  \mathrm{}\, \textit{}

\\

\\
&=\displaystyle{\int_{\mathcal G^0}h(u)\int_{\mathcal G_1\times \mathcal G_2}(f\overline{f'})({g}_1^{-1}{g}_2)\delta_1(g_1^{-1})(\delta\delta_2^{-1})({g }_2) d\lambda_2^{u}(g_2)d\lambda_1^{u}(g_1)d\nu(u)} & \mathrm{}\, \textit{}

\end{array}$$
 the lemma follows.
  \end{dm}
\vskip 0.6cm

\subsubsection{\bf Notations}
\label{represent}
For any $i,j\in \{s,r,m,\hat m\}$,   let  us denote  $\mu^2_{i,j} = \mu _{i_{\mathcal G}}\times_{j_{\mathcal G}} \mu$ and $\mathcal  G_{i,j}= \{(g,g') \in \mathcal G/ i(g) = j(g')\}$,   therefore   $L^2(\mathcal G,\mu)\underset{\nu}{_{i_{\mathcal G}}\otimes_{j_{\mathcal G}}}L^2(\mathcal G ,\mu)$  is isomorphic to $L^2(\mathcal G^2_{i,j},\mu^2_{i,j})$.\vskip 0.6cm

\subsubsection{\bf Lemma}
\label{metr}
{\it    For any $f \in \mathcal K(\mathcal G\times \mathcal G)$ one has: 
$$ \mu^2_{m,r}(f) = \int _{\mathcal G^0}\int_{\mathcal G^2_{m,r} }f(g,g')d\lambda^{m(g)}(g')d\lambda^u(g)d\nu(u) $$
$$ \mu^2_{s,m}(f) = \int _{\mathcal G^0}\int_{\mathcal G^2_{s,m} }\delta(g^{-1}) f(g^{-1},{g'} )d\lambda^{m(g')}(g)d\lambda^u(g')d\nu(u),  $$ }
\vskip 0.2cm
\begin{dm}  This is    easy computations.
\end{dm}
\vskip 0.8cm

Let's define an important   pseudo multiplicative unitary, which  generalizes at the same time,   the multiplicative unitary   of \cite{BSV1} 3.2 and the multiplicative partial isometry $I_{\mathcal H,\mathcal K}$ of \cite{Val2} definition 4.1.5.

\vskip 1cm

\subsubsection{\bf Proposition}
\label{trombone}
{\it Let $W_{\mathcal G_1,\mathcal G_2} : L^2(\mathcal G^2_{s,m},\mu^2_{s,m}) \to L^2(\mathcal G^2_{m,r},\mu^2_{m,r})$ be the operator defined for any $\xi \in  L^2(\mathcal G^2_{s,m},\mu^2_{s,m})$ and $\mu^2_{m,r}$-almost any $(x,y)$ in $G^2_{m,r}$ by:
$$  W_{\mathcal G_1,\mathcal G_2}\xi(x,y) = D(x,y)^{\frac{1}{2}}\xi(\theta(x,y))$$
where $\theta(x,y) = (xp_1(p_2(x)^{-1}y), p_2(x)^{-1}y)$ and $D(x,y)$ is  the Radon Nikodym derivative $\frac{d\mu^2_{s,m}\circ \theta}{d\mu^2_{m,r}}$ .  Up to the identification of \ref{represent}, $W_{\mathcal G_1,\mathcal G_2} $ is a pseudo multiplicative unitary over the basis $(L^\infty(\mathcal G^0,\nu),m_{\mathcal G},s_{\mathcal G},r_{\mathcal G})$.}
\vskip 0.7cm
\begin{dm} Mimicking  \cite{BSV1} 3.2, let us consider two maps,   $\omega_1$ and $\omega_2$, defined on $G^2_{m,r} \cap \mathcal G_2\mathcal G_1\times \mathcal G$ and on  $G^2_{s,r} \cap  \mathcal G\ \times \mathcal G_2\mathcal G_1$ respectively by: $\omega_1(x,y)= (x,p_2(x)^{-1}y)$  and $\omega_2(x,y)  =(xp_1(y),y)$). Obviously, one has: $Im \omega_1 \subset \mathcal G^2_{s,r}$ (resp. $Im \omega_2 \subset \mathcal G^2_{s,m}$). Due to   lemma \ref{metr},  for any $f \in \mathcal K(\mathcal G \times \mathcal G)$, one has:
$$ \mu^2_{m,r}(f\circ \omega_1) =  \mu^2_{s,r}(f), \hskip 1cm \mu^2_{s,r}(f\circ \omega_2) =  \mu^2_{s,m}(\Delta f) $$where $\Delta(x,y)= \delta(p_1(y))$. So one can define two unitaries  $W_1: L^2(\mathcal G^2_{s,r},\mu^2_{s,r}) \to L^2(\mathcal G^2_{m,r},\mu^2_{m,r})$ and $W_2:L^2(\mathcal G^2_{s,m},\mu^2_{s,m} )\to L^2(\mathcal G^2_{s,r},\mu^2_{s,r})$ such that $W_1\xi = \xi \circ \omega_1$ and $W_2\eta = \Delta^{-1/2}\eta\circ \omega_2$. Their composition $ W_{\mathcal G_1,\mathcal G_2} = W_1W_2$   is  defined for any $\xi \in  L^2(\mathcal G^2_{s,m},\mu^2_{s,m})$ and $\mu^2_{m,r}$-almost any $(x,y)$ in $G^2_{m,r}$ by:
$$  W_{\mathcal G_1,\mathcal G_2}\xi(x,y) = D(x,y)^{\frac{1}{2}}\xi(\theta(x,y))$$ where $D$ is  the Radon Nikodym derivative $\frac{d\mu^2_{s,m}\circ \theta}{d\mu^2_{m,r}}$ and $\theta (x,y) = \omega _2 \circ\omega_1(x,y) = (xp_1(p_2(x)^{-1}y), p_2(x)^{-1}y) $. $W_{\mathcal G_1,\mathcal G_2}$ is obviously a unitary and the fact that this is a pseudo multiplicative unitary is essentially proposition 4.1.6 in \cite{Val2}.
\end{dm}

\vskip 0.2cm

\subsubsection{\bf Remark}
Using the identification of \ref{sarcofobe},   $W_{\mathcal G_1,\mathcal G_2}$ is also a unitary:  
\begin{align*}
[L^2(\mathcal G_2,\mu_2)&\underset{\nu}{_{s_2}\otimes_{r_1}}L^2(\mathcal G_1,\mu_1)]\underset{\nu}{_{s_1}\otimes_{m}}[L^2(\mathcal G_2,\mu_2)\underset{\nu}{_{s_2}\otimes_{r_1}}L^2(\mathcal G_1,\mu_1)] 
\\
&\to [L^2(\mathcal G_2,\mu_2)\underset{\nu}{_{s_2}\otimes_{r_1}}L^2(\mathcal G_1,\mu_1)]\underset{\nu}{_{m }\otimes_{r_2}}[L^2(\mathcal G_2,\mu_2)\underset{\nu}{_{s_2}\otimes_{r_1}}L^2(\mathcal G_1,\mu_1)] 
\end{align*}

\vskip 0.6cm

\subsubsection{\bf{Notation}}  
\label{vraicourage}
{\it Due to \ref{bienfait} and \ref{sarcofobe}, one can consider the fiber product
$$\mathfrak a \star \mathfrak a : L^\infty(\mathcal G_2,\mu_2)\underset{\nu}{_{s_2}\star_{r_2}}L^\infty(\mathcal G_2,\mu_2) \to L^\infty(\mathcal G ,\mu )\underset{\nu}{_{s }\star_{ m}}L^\infty(\mathcal G ,\mu )$$}

\subsubsection{\bf{Lemma }}  
\label{blair}
{\it For $\mu^2_{s,m}$-almost any $(g,g') \in \mathcal G^2_{s,m}$, one has: 
$$D(g,g') = \delta^{-1}(p_1(p_2(g)^{-1}g'))$$}
\label{tillend}
\vskip 0.2cm
\begin{dm}
For  $\mu^2_{s,m}$-almost any $(g,g') \in \mathcal G^2_{s,m}$, and any $\xi \in  L^2(\mathcal G^2_{s,m},\mu^2_{s,m})$ one has:
\begin{align*}
W\xi(g,g')& = W_1W_2(g,g') = W_2\xi(g,p_2(g)^{-1}g') \\
&= \Delta^{-\frac{1}{2}}(g,p_2(g)^{-1}g')\xi(gp_1(p_2(g)^{-1}g'), p_2(g)^{-1}g' )\\
&= \delta^{-\frac{1}{2}}(p_1(p_2(g)^{-1}g'))\xi(gp_1(p_2(g)^{-1}g'), p_2(g)^{-1}g' )
\end{align*}
The lemma follows.
\end{dm}

\vskip 1.1cm

\subsubsection{\bf{Proposition}}  
\label{dur}
 {\it   The von Neumann algebra $M$ (resp.$\hat M$) generated by the left (resp.right) leg  of  $W_{\mathcal G_1,\mathcal G_2} $ is isomorphic to the crossed product $L^\infty(\mathcal G_2,\mu_2) \rtimes_{\mathfrak a} \mathfrak G(\mathcal G_1)$ (resp.$L^\infty(\mathcal G_1,\mu_1) \rtimes_{\hat{\mathfrak a}} \mathfrak G(\mathcal G_2)$).  }
\vskip 0.3cm
\begin{dm} 
For any $f,h,\eta, \eta'  \in \mathcal K( \mathcal G)$, in $L^2(\mathcal G,\mu)$ one has: 
 \begin{align*}
 ((i \star &\omega_{ f,h})(W_{\mathcal G_1,\mathcal G_2})\eta,\eta') 
 =    \int_{\mathcal G^2_{m,r}} W_{\mathcal G_1,\mathcal G_2}(\eta \underset{\nu}{_{s_2}\otimes_{r_1}} f )(g,g')\overline{\eta' (g) h(g')}d\mu^2_{m,r}(g,g')\\
 &= \int_{\mathcal G^2_{m,r}}D(g,g')^{\frac{1}{2}}\eta(gp_1(p_2(g)^{-1}g'))f(p_2(g)^{-1}g')\overline{\eta' (g) h(g')}d\mu^2_{m,r}(g,g')\\
 &= \int_{\mathcal G^0} \int_{\mathcal G^2_{m,r}}D(g,g')^{\frac{1}{2}}\eta(gp_1(p_2(g)^{-1}g'))f(p_2(g)^{-1}g')\overline{\eta' (g) h(g')} d\lambda^{m(g)}(g')d\lambda^u(g)d\nu(u) \\
 &= \int_{\mathcal G}(\int_{\mathcal G}D(g,g')^{\frac{1}{2}}\eta(gp_1(p_2(g)^{-1}g'))f(p_2(g)^{-1}g')\overline{h(g')}d\lambda^{m(g)}(g'))\overline{\eta' (g)} d\mu(g) 
 \end{align*}

Let's change  of variable: $g' \mapsto p_2(g)^{-1}g'$:

\begin{align*}
 ((i \star &\omega_{ f,h})(W_{\mathcal G_1,\mathcal G_2})\eta,\eta') =  \\
&= \int_{\mathcal G}(\int_{\mathcal G} D(g,p_2(g)g')^{\frac{1}{2}}\eta (gp_1(g')f(g')\overline h(p_2(g)g')d\lambda^{s(g)}(g'))\overline{\eta' (g)} d\mu(g)   
\end{align*}

which gives, using Proposition \ref{use} and lemma \ref{blair} , that for $\mu$-almost any $g'\in \mathcal G$, one has:

\begin{align*}
 &(i \star \omega_{ f,h})(W_{\mathcal G_1,\mathcal G_2})\eta(g)=   \\
&=  \int_{\mathcal G_1}(\int_{\mathcal G_2} D(g,p_2(g)g'_1g'_2)^{\frac{1}{2}}\eta (gg'_1)f(g'_1g'_2)\overline h(p_2(g)g'_1g'_2)\frac{\delta(g'_2)}{\delta_2(g'_2)}) d\lambda_2^{s(g'_1)}(g'_2)d\lambda^{s(g)}_1(g'_1)  \\
&=  \int_{\mathcal G_1}(\int_{\mathcal G_2} \delta(g'_1)^{-\frac{1}{2}}\eta (gg'_1)f(g'_1g'_2)\overline h(p_2(g)g'_1g'_2)\frac{\delta(g'_2)}{\delta_2(g'_2)}) d\lambda_2^{s(g'_1)}(g'_2)d\lambda^{s(g)}_1(g'_1) 
\end{align*}

Let us   note: $\Theta(g,g'_1,g'_2) = f(g'_1g'_2)\overline h(gg'_1g'_2)\frac{\delta}{\delta_2}(g'_2)     $,  then:
 
\begin{align*}
 (i \star \omega_{ f,h})(W_{\mathcal G_1,\mathcal G_2})\eta(g)  
&=  \int_{\mathcal G_1}(\int_{\mathcal G_2}\Theta(p_2(g),g'_1,g'_2)  d\lambda_2^{s(g'_1)}(g'_2))\delta(g'_1)^{-\frac{1}{2}} \eta (gg'_1)d\lambda^{s(g)}_1(g'_1) \\
&= \mathfrak R(F_{f,h})\eta(g)
\end{align*}

where $F_{f,h}(g_2,g_1) = \int_{\mathcal G_2}\theta(g_2,g_1,g'_2)d\lambda_2^{s(g_1)}(g'_2)$

This implies  that  the left leg of $W_{\mathcal G_1,\mathcal G_2} $ generates the crossed product $L^\infty(\mathcal G_2,\mu_2) \rtimes_{\mathfrak a} \mathfrak G(\mathcal G_1)$, analogue computations give that its right leg generates $L^\infty(\mathcal G_1,\mu_1) \rtimes_{\hat{\mathfrak a}} \mathfrak G(\mathcal G_2)$.
\end{dm}
\vskip 0.5cm

\vskip 0.3cm
 Using Proposition \ref{dur}, we shall identify the left (resp.right) leg of $W_{\mathcal G_1,\mathcal G_2} $ with  crossed products.
\vskip 0.5cm

\subsubsection{\bf{Corollary }}
\label{logaplouf}  
{\it Thanks to the existence of $W_{\mathcal G_1,\mathcal G_2} $, one can define two Hopf bimodule structures. We shall note them   $(L^\infty(\mathcal G^0,\nu),L^\infty(\mathcal G_2,\mu_2) \rtimes_{\mathfrak a} \mathfrak G(\mathcal G_1),m,s,\Gamma   )$ for the left leg,  and $(L^\infty(\mathcal G^0,\nu),L^\infty(\mathcal G_1,\mu_1) \rtimes_{\hat{\mathfrak a}} \mathfrak G(\mathcal G_2) ,r,m,\hat \Gamma $) for the right one.}
\vskip 0.1cm
\begin{dm}
This is a consequence of remark \ref{echocard}
\end{dm}

\vskip 1cm

\section{The quantum groupoid structures associated with a matched pair}
In this chapter we shall describe in details the   structures found in the previous one. We shall  complete them to obtain measured quantum groupoids structures. In order to simplify notations and using  \ref{sarcofobe}, we  can suppose that  $L^\infty(\mathcal G_2,\mu_2) \rtimes_{\mathfrak a} \mathfrak G(\mathcal G_1)$ is acting on   $L^2(G,\mu)$  and that it is generated by products $\mathfrak a(f)(1 {_{s_2}\otimes}_{r_1}\rho( h))$, where, for any   $f \in L^\infty(\mathcal G_2,\mu_2)$ and $\mu$-almost  any $g\in \mathcal G$,  one has  $\mathfrak a(f)(g) = f(p_2(g))$ and  for any $\xi \in L^2(\mathcal G,\mu)$, $h \in  \mathcal K (\mathcal G_1)$  : $(1 {_{s_2}\otimes}_{r_1}\rho( h))\xi(g) = \int_{\mathcal G_1}h(g_1)\xi(gg_1)d\lambda_1^{s(g)}(g_1)$.

\subsection{The coproduct } 
\subsubsection{\bf Lemma }
\label{couscous1}
{\it One has: $\Gamma \circ \mathfrak a= (\mathfrak a  _{s_2}\star_{r_2}  \mathfrak a)\Gamma_{\mathcal G_2}$, hence for any $f \in L^\infty(\mathcal G_2,\mu_2)$   and $\mu^2_{s,m}$-almost any $(g,g') \in \mathcal G^2_{s,m}$  , one gets: 
$$\Gamma(\mathfrak a(f)) (g,g') =f(p_2(g)p_2(g'))    $$  }
\newline
\begin{dm}
As $W_{\mathcal G_1,\mathcal G_2}  $ is a unitary, it is easy to see that for any  $\eta \in L^2(\mathcal G^2_{m,r},\mu^2_{m,r})$ and $\mu^2_{s,m}$-almost any $(g,g') \in \mathcal G^2_{s,m}$, one has: 
$$W_{\mathcal G_1,\mathcal G_2}^* \eta(g,g') = \delta(p_1(g'))^{\frac{1}{2}}\eta(gp_1(g')^{-1},p_2(gp_1(g')^{-1})g')$$
As $D'(g,g') = \delta(p_1(g'))^{\frac{1}{2}}$ is a  density, then $D'(g,g')D((gp_1(g')^{-1},p_2(gp_1(g')^{-1})g') = 1$,  hence
for any $f \in L^\infty(\mathcal G_2,\mu_2)$, any $\xi \in L^2(\mathcal G^2_{s,m},\mu^2_{s,m})$ and $\mu^2_{s,m}$-almost any $(g,g') \in \mathcal G^2_{s,m}$, one has: 
\begin{align*}
\Gamma(\mathfrak a(f)) \xi(g,g') 
&=W_{\mathcal G_1,\mathcal G_2}^* (1{ _m\otimes_r}\mathfrak a(f)) W_{\mathcal G_1,\mathcal G_2}\xi(g,g') \\
&= D'(g,g')^{-\frac{1}{2}}(1{ _m\otimes_r}\mathfrak a(f))W_{\mathcal G_1,\mathcal G_2}\xi(gp_1(g')^{-1},p_2(gp_1(g')^{-1})g')\\
&=\mathfrak a(f)(p_2(gp_1(g')^{-1})g')\xi(g,g') = f(p_2(p_2(gp_1(g')^{-1})g'))\xi(g,g') \\
&= f(p_2(gp_1(g')^{-1}g'))\xi(g,g') = f(p_2(gp_2(g'))\xi(g,g') \\
&= f(p_2(g)p_2(g'))\xi(g,g')
\end{align*} 
\end{dm}
\vskip 0.5cm

A good description of  $\Gamma( 1 {_{s_2}\otimes}_{r_1} \mathcal R(\mathcal G_1))$ is   given by  an integral.

\vskip 1cm

\subsubsection{\bf Proposition }
\label{couscous 2}
{\it   Let   $h$ (resp.$f$)  be any element in $ \mathcal K (\mathcal G_1)$  (resp.$ \mathcal K(\mathcal G_2)$), then: 

i) for   all  $\xi  \in \mathcal K(\mathcal G^2_{s,m})$ and $\mu^2_{s,m}$-almost any $(g,g') \in \mathcal G^2_{s,m}$, one has:
$$   \Gamma(\mathfrak a(f)(1 {_{s_2}\otimes}_{r_1} \rho( h))\xi  (g,g') =  f(p_2(g)p_2(g'))   \int_{\mathcal G_1}h(g_1) \xi(gp_1(p_2(g')g_1),g'g_1)d\lambda_1^{s(g')}(g_1)$$
ii) with the notations of \ref{dublin}, for any $\phi \in \mathcal K(\mathcal G)$:
$$(\omega_{ \phi} \underset{\nu}{_{s }\star_{m}}i)( \Gamma(\mathfrak a(f)(1 {_{s_2}\otimes}_{r_1} \rho( h)))) = \mathfrak R(\Psi_{f,h})$$ 
where $\Psi_{f,h} \in \mathcal K(G_2   {_{s_2}\times_{r_1}} \mathcal G_1)$ is defined for any $(g_2,g_1) \in G_2   {_{s_2}\times_{r_1}} \mathcal G_1$ by:
$$\Psi_{f,h}(g_2,g_1) = h(g_1)\int _{\mathcal G^{r(g_2)}}\delta(g^{-1}) f(p_2(g^{-1})g_2)\phi(g^{-1}p_1(g_2g_1))\overline{\phi(g^{-1})}d\lambda^{r(g_2)}(g)$$}
 \vskip 0.5cm
 \begin{dm}
i) For any   $h \in  \mathcal K (\mathcal G_1)$, any $\xi \in \mathcal K(\mathcal G^2_{s,m})$ and $\mu^2_{s,m}$-almost any $(g,g') \in \mathcal G^2_{s,m}$, one has: 
\begin{align*}
\Gamma(1 {_{s_2}\otimes}_{r_1}&\rho( h))\xi  (g,g') 
= W_{\mathcal G_1,\mathcal G_2}^*(1 {_m\otimes}_{r}(1 {_{s_2}\otimes}_{r_1}\rho( h)))W_{\mathcal G_1,\mathcal G_2} \xi  (g,g')  \\
&= D'(g,g')^{-\frac{1}{2}}(1 {_m\otimes}_{r}(1 {_{s_2}\otimes}_{r_1}\rho( h)))W_{\mathcal G_1,\mathcal G_2}\xi  (gp_1(g')^{-1},p_2(gp_1(g')^{-1})g')\\
&=  \int_{\mathcal G_1}h(g_1) \xi(gp_1(p_2(g')g_1),g'g_1)d\lambda_1^{s(g')}(g_1)
\end{align*}

i) follows immediately.

ii) For any $f \in \mathcal K(\mathcal G_2),\phi,\phi' \in \mathcal K(\mathcal G)$ and  $\mu^2_{s,m}$-almost any  $(g,g') \in \mathcal G^2_{s,m}$,   let us note: 
$$X_{\phi,\phi'}^f(g,g') =  f(p_2(g)p_2(g')) \int_{\mathcal G^1}  h(g_1)  \phi(g p_1(p_2(g')g_1))\phi'(g'g_1)d\lambda_1^{s(g')}(g_1).$$
Due to i) and \ref{couscous1}, one has:
  
 $$\begin{array}{cll}
 &(\omega_\phi \underset{\nu}{_{s }\otimes_{m}}\omega_{ \phi'}) ( \Gamma(\mathfrak a(f)(1 {_{s_2}\otimes}_{r_1} \rho( h))))=
\\

\\
 &= \int_{\mathcal G^2_{s,m}}  \int_{\mathcal G_1}h(g_1) f(p_2(g)p_2(g'))\phi(gp_1(p_2(g')g_1))\phi'(g'g_1)d\lambda_1^{s(g')}(g_1)\overline{\phi(g)\phi'(g')}d\mu^2_{s,m} 
 \\

\\
 &= \int_{\mathcal G^2_{s,m}} X_{\phi,\phi'}^f(g,g') \overline{\phi(g)\phi'(g')}d\mu^2_{s,m}(g,g')
 \\

\\
  &= \int_{\mathcal G}\Big( \int_{\mathcal G^ {m(g')}} \delta(g^{-1})X_{\phi,\phi'}^f(g^{-1},g') \overline{\phi(g^{-1})}d\lambda^{ m (g')}(g)\Big)\overline {\phi'(g')}d\mu(g')
 
  \quad \mathrm{by}\, \textit{\ref{metr}}\\

\end{array}$$
\vskip 1.1cm

If    $\Psi_{f,h} \in \mathcal K(G_2   {_{s_2}\times_{r_1}} \mathcal G_1)$ is defined for any $(g_2,g_1) \in G_2   {_{s_2}\times_{r_1}} \mathcal G_1$ by:
$$\Psi_{f,h}(g_2,g_1) = h(g_1)\int _{\mathcal G^{r(g_2)}}\delta(g^{-1}) f(p_2(g^{-1})g_2)\phi(g^{-1}p_1(g_2g_1))\overline{\phi(g^{-1})}d\lambda^{r(g_2)}(g)$$

we can write:

$$\begin{array}{cll}
& \int_{\mathcal G^ {m(g')}} \delta(g^{-1})X^f_{\phi,\phi'}(g^{-1},g') \overline{\phi(g^{-1})}d\lambda^{ m (g')}(g)
=
 \\

\\
&=\int_{\mathcal G_1}h(g_1) \int_{\mathcal G^ {m(g')}}f(p_2(g^{-1})p_2(g'))  \phi(g^{-1} p_1(p_2(g')g_1))\overline{\phi(g^{-1})}d\lambda^{ m (g')}(g)\phi'(g'g_1)d\lambda_1^{s(g')} 

 \\

\\
&=\int_{\mathcal G_1}\Psi_{f,h}(p_2(g'),g_1)\phi'(g'g_1)d\lambda_1^{s(g')}(g_1) \  \  \   \quad \mathrm{as \  \ m(g') = r(p_2(g'))} 
\end{array}$$

\vskip 2cm
hence one has:

$$\begin{array}{cll}
 &\omega_{ \phi'} \Big ((\omega_\phi \underset{\nu}{_{s }\star_{m}}i) ( \Gamma(1 {_{s_2}\otimes}_{r_1} \rho( h)))\Big ) = (\omega_\phi \underset{\nu}{_{s }\star_{m}}\omega_{ \phi'}) ( \Gamma(1 {_{s_2}\otimes}_{r_1} \rho( h)))
 \\

 \\
 &= \int_{\mathcal G}\Big(\int_{\mathcal G_1}\Psi_{f,h}(p_2(g'),g_1)\phi'(g'g_1)d\lambda_1^{s(g')}(g_1) \Big)\overline {\phi'(g')}d\mu(g')
 
 \\

\\
  &= \int_{\mathcal G} \Big( \mathfrak R(\Psi_{f,h})\phi'(g')\Big)\overline {\phi'(g')}d\mu(g') = \omega_{\phi'}(\mathfrak R(\Psi_{f,h}))
\end{array}$$
ii) follows
\end{dm}

\subsubsection{\bf{Remark}}
\label{prg}
The formulas of lemma \ref{couscous 2} generalize the ones obtained by Stefaan Vaes in \cite{Vas1} 4.20 (and maybe elsewhere). 

\vskip 2cm
\subsection{The co-involution  } In this paragraph, we define a co-involution on $ L^\infty(\mathcal G_2,\mu_2) \rtimes_{\mathfrak a} \mathfrak G(\mathcal G_1) $.
\vskip 0.2cm

\subsubsection{\bf{Lemma}}
\label{georgie}
{\it Using notations of Remark \ref{rollex}, let $\phi :    \mathcal G \to \mathcal G$   (resp. $\hat \phi : \mathcal G \to \mathcal G$) be the   function defined by  the formula: $\phi(g) =  {p_1}(g )^{-1} p'_2(g) $  (resp. $\hat \phi (g) =  \phi(g)^{-1}$), then one has:

i) $\hat \phi(g) = {p'_2(g )}^{-1} p_1(g) =p'_1(g){p_2(g)}^{-1} $

ii) $\phi^2 = \hat \phi^2 = id_{\mathcal G}$

iii) $ s \circ \phi = \hat m, m \circ \phi = r, m \circ \hat \phi = s,  r  \circ  \phi = m $}
\vskip 0.2cm
\begin{dm} 
For  any $ g  \in \mathcal G_1\mathcal G_2 \cap \mathcal G_2\mathcal G_1$, one has $p_1(g)p_2(g) = p'_2(g)p'_1(g)$ so i) is true; using obvious notations, let us write: $g= g_1g_2 = g'_2g'_1$; one has: $\phi(\phi(g)) = \phi(g_1^{-1}g'_2) = \phi(g_2{g'_1}^{-1})=  g_1g_2  = g$, this gives ii), the assertion iii) is obvious.
\end{dm}

\vskip 1.1cm
Using lemma \ref{georgie}, one can define a map  $\hat \phi {_{s}\times}_{m} \phi : \mathcal G^2_{s,m} \to  \mathcal G^2_{m,r}$ (resp $\hat \phi {_m\times}_{r} \phi : \mathcal G^2_{m,r} \to  \mathcal G^2_{s,m}$) such that for almost any $(g,g') \in \mathcal G^2_{s,m}$ (resp.$\mathcal G^2_{m,r}$), one has:  $(\hat \phi {_s\times}_{m} \phi)(g,g') = (\hat \phi(g), \phi(g'))$  (resp. $(\hat \phi {_m\times}_{r} \phi)(g,g') = (\hat \phi(g), \phi(g'))$.

\subsubsection{\bf{Lemma}}
\label{ossetie}
{\it Using the notations of \ref{georgie} and \ref{trombone}, one has:
$$  \theta(\hat \phi {_s\times}_{m}  \phi) = (\hat \phi {_m\times}_{r} \phi)\theta^{-1}.$$ }
\vskip 0.03cm
\begin{dm} For almost  any $(g,g') \in \mathcal G^2_{m,r}$, due to \ref{georgie} i), one has   $p_2(\hat \theta(g)) =  p_2(g)^{-1}$, so :
\begin{align*}
\theta(\hat \phi {_s\times}_{m}  \phi)(g,g') 
&= \theta(\hat \phi(g),\phi(g'))   =  ( \hat \phi(g)p_1(p_2(\hat \phi(g))^{-1}\phi(g')),p_2(\hat \phi(g))^{-1}\phi(g')) \\
&= ( {p'_2(g )}^{-1} p_1(g)p_1(p_2(g)\phi(g')),p_2(g)\phi(g')) \\
&= ( {p'_2(g )}^{-1}  p_1(g{p_1}(g' )^{-1} p'_2(g')),p_2(g){p_1}(g' )^{-1} p'_2(g') ) \\
&= ( {p'_2(g )}^{-1}  p_1(g{p_1}(g' )^{-1} ),p_2(g){p_1}(g' )^{-1} p'_2(g') )
\end{align*}

Also one has:
\begin{align*}
 (\hat \phi {_m\times}_{r} \phi)\theta^{-1}(g,g')
 &=  (\hat \phi {_m\times}_{r} \phi)(gp_1(g')^{-1},p_2(gp_1(g')^{-1})g')
\end{align*}

Let us define:
$$ X= gp_1(g')^{-1}, \  \ Y= p_2(gp_1(g')^{-1})g'$$
this gives that:
\begin{align*}
 (\hat \phi {_m\times}_{r} \phi)\theta^{-1}(g,g') 
 &= (p'_2(X)^{-1}p_1(X),p_1(Y)^{-1}p'_2(Y)) \\
 &= (p'_2(X)^{-1}p_1(X),p_2(Y)p'_1(Y)^{-1})
\end{align*}
\vskip 0.6cm
as $p'_2(X) = p'_2(g)$, one deduces that: 

$$p'_2(X)^{-1}p_1(X) = {p'_2(g )}^{-1}  p_1(g{p_1}(g' )^{-1} )$$
\vskip 0.6cm
Also we can write that:
\begin{align*}
Y[p_2(g)p_2(g')]^{-1}
&= p_2(g{p_1}(g' )^{-1} )p_1(g')p_2(g')[p_2(g)p_2(g')]^{-1}= p_2(g{p_1}(g' )^{-1} )p_1(g')p_2(g)^{-1}\\
&= p_2(p_2(g){p_1}(g' )^{-1} )[p_2(g)p_1(g')^{-1}]^{-1}\\
&= [ p_1(p_2(g){p_1}(g' )^{-1} )]^{-1}
\end{align*}

this implies that  $ p_2(Y) = p_2(g)p_2(g')$, 
and as $Y= p_2(gp_1(g')^{-1})g'$,  therefore  $p'_1(Y) = p'_1(g')$, which gives that:
 \vskip 0.6cm
 
 $$  p_2(Y)p'_1(Y)^{-1} =  p_2(g)p_2(g')p'_1(g')^{-1}$$
 
 finally:
 
 \begin{align*}
  (\hat \phi {_m\times}_{r} \phi)\theta^{-1}(g,g') 
  &= (p'_2(X)^{-1}p_1(X),p_2(Y)p'_1(Y)^{-1}) \\
  &= ({p'_2(g )}^{-1}  p_1(g{p_1}(g' )^{-1}), p_2(g)p_2(g')p'_1(g')^{-1})\\
  &=\theta(\hat \phi {_s\times}_{m}  \phi)(g,g') 
 \end{align*}
 
which gives the lemma.
\end{dm}

\vskip 0.4cm

\subsubsection{\bf{Proposition}}
\label{majid}
{\it i) Let $J,\hat J : L^2(\mathcal G) \to  L^2(\mathcal G)$  be defined for any $\xi \in L^2(\mathcal G,\mu)$ and $\mu$-almost any $g \in \mathcal G$ by :  $$  J\xi(g) = \overline{ \xi(  \phi(g))}  (\frac{d\mu \circ \phi}{d\mu})^{\frac{1}{2}}(g),\hskip 0.5cm \hat J\xi(g )  = \overline{ \xi(\hat \phi(g) )} (\frac{d\mu \circ \hat  \phi}{d\mu})^{\frac{1}{2}}(g)$$ 
 
 then $J$ and $\hat J$ are antilinear involutive isometries.
 
 ii) For all $f \in \mathcal K(\mathcal G^0)$, one has: $\hat J s(f) = m(\overline f) \hat J$ and $J m(f) = r(\overline f)J$, hence one can define $\hat J  {_{s}\otimes}_{m}J : L^2(\mathcal G^2_{s,m},\mu^2_{s,m}) \to L^2(\mathcal G^2_{m,r},\mu^2_{m,r})$, and it's inverse $\hat J  {_{m}\otimes}_{r}J : L^2(\mathcal G^2_{m,r},\mu^2_{m,r}) \to L^2(\mathcal G^2_{s,m},\mu^2_{s,m})$, which  verify:
 
$$(\hat J  {_{m}\otimes}_{r}J) W_{\mathcal G_1,\mathcal G_2}  =  W_{\mathcal G_1,\mathcal G_2}^*(\hat J  {_{s}\otimes}_{m}J) $$}
\vskip 0.03cm
\begin{dm}
This is a straightforward consequence of lemma \ref{ossetie}.\end{dm}

\vskip 1cm
\subsubsection{\bf{Proposition}}
\label{majid2}
{\it Let $R$ be the application defined for any  $x \in L^\infty(\mathcal G_2,\mu_2) \rtimes_{\mathfrak a} \mathfrak G(\mathcal G_1)$ by: $R(x) =  \hat Jx^*\hat J$, then $R$ is a co-involution of  $(L^\infty(\mathcal G^0,\nu),L^\infty(\mathcal G_2,\mu_2) \rtimes_{\mathfrak a} \mathfrak G(\mathcal G_1),m,s,\Gamma)$, more precisely, for any $f \in \mathcal K(\mathcal G_2)$, if we define $f^{-1}$   by $f^{-1}(g_2) = f({g_2}^{-1})$ one has: $R(\mathfrak a(f))= \mathfrak a(f^{-1})$.}
\vskip 0.3cm
\begin{dm}
Obviously $Ad(\hat J)$ is an involutive $*$-antiautomorphism of $\mathcal L(L^2(\mathcal G))$. Due to proposition \ref{majid}, for any $h,k \in \mathcal K(\mathcal G)$, one has: 
\begin{align*}
R( (i \star \omega_{h,k})(W_{\mathcal G_1,\mathcal G_2})) 
&=\hat J (i \star \omega_{h,k})(W_{\mathcal G_1,\mathcal G_2})^*\hat J    = \hat J (i \star \omega_{k,h})(W_{\mathcal G_1,\mathcal G_2}^*) \hat J \\
&=  (i \star \omega_{Jk,Jh})(\hat J  {_{m}\otimes}_{r}J)W_{\mathcal G_1,\mathcal G_2}^*(\hat J  {_{m}\otimes}_{r}J) )  =  (i \star \omega_{Jk,Jh})(W_{\mathcal G_1,\mathcal G_2} ) 
\end{align*}
so, by restriction, $R$ is an involutive $*$-antiautomorphism of $ L^\infty(\mathcal G_2,\mu_2) \rtimes_{\mathfrak a} \mathfrak G(\mathcal G_1)$.

From proposition \ref{majid}, one obtains   that $R \circ m = s$.

For all $f \in \mathcal K(\mathcal G_2)$,  all $\xi \in L^2(\mathcal G,\mu)$ and $\mu$-almost any $g \in \mathcal G_1\mathcal G_2 \cap \mathcal G_2\mathcal G_1$, with $g= g_1g_2 = g'_1g'_2$, one has: 
\begin{align*}
R(\mathfrak a(f))\xi(g) 
&= \hat J\mathfrak a^*(f)\hat J\xi(g) = \overline{\mathfrak a (\overline f)(\hat \phi(g)) \xi( g)} = f(p_2({g'_2}^{-1}g_1))\xi(g) = f(p_2(g'_1g_2^{-1})\xi(g) \\
&= f(g_2^{-1})\xi(g) = \mathfrak a(f^{-1})(g)\xi(g)
\end{align*} 

So $R(\mathfrak a(f)) = \mathfrak a(f^{-1})$.

Due to Proposition 3.7 of \cite{E3}, for any $k,k_1,f_2,h_1,h_2 \in \mathcal K(\mathcal G)$, one has:
$$( \Gamma((i \star \omega_{h,k})(W_{\mathcal G_1,\mathcal G_2}))(h_1\underset{\nu}{_{s}\otimes_{m}}k_1),h_2\underset{\nu}{_{s}\otimes_{m}}k_2)= ((\omega_{h_1,h_2}\star i)(W_{\mathcal G_1,\mathcal G_2})(\omega_{k_1,k_2}\star i)(W_{\mathcal G_1,\mathcal G_2})h,k)$$ 
Hence, in the one hand:

\begin{align*}
(\varsigma_{m,s} \Gamma(R((i \star \omega_{h,k})&(W_{\mathcal G_1,\mathcal G_2})))(h_1\underset{\nu}{_{m}\otimes_{s}}k_1),h_2\underset{\nu}{_{m}\otimes_{s}}k_2)\\
&= ( \Gamma( (i \star \omega_{Jk,Jh}) (W_{\mathcal G_1,\mathcal G_2})))(k_1\underset{\nu}{_{s}\otimes_{m}}h_1),k_2\underset{\nu}{_{s}\otimes_{m}}h_2)\\
&= ((\omega_{k_1,k_2}\star i)(W_{\mathcal G_1,\mathcal G_2})(\omega_{h_1,h_2}\star i)(W_{\mathcal G_1,\mathcal G_2})Jk,Jh)  
\end{align*}
 On the other hand, using proposition \ref{majid} ii), one has:
 
 \begin{align*}
 ((R  {_{m}\star_{s}} R)\Gamma((i \star \omega_{h,k})&(W_{\mathcal G_1,\mathcal G_2}))(h_1\underset{\nu}{_{m}\otimes_{s}}k_1),h_2\underset{\nu}{_{m}\otimes_{s}}k_2) \\
 &=(  \hat Jh_2\underset{\nu}{_{s}\otimes_{m}}\hat Jk_2,\Gamma((i \star \omega_{h,k})(W_{\mathcal G_1,\mathcal G_2})^*)(\hat Jh_1\underset{\nu}{_{s}\otimes_{m}}\hat Jk_1))\\
 &= ( \Gamma((i \star \omega_{h,k})(W_{\mathcal G_1,\mathcal G_2}) )(\hat Jh_2\underset{\nu}{_{s}\otimes_{m}}\hat Jk_2),\hat Jh_1\underset{\nu}{_{s}\otimes_{m}}\hat Jk_1)\\
 &= ((\omega_{\hat Jh_2,\hat Jh_1}\star i)(W_{\mathcal G_1,\mathcal G_2})(\omega_{\hat Jk_2, \hat Jk_1}\star i)(W_{\mathcal G_1,\mathcal G_2})h,k) \\
  &= (J(\omega_{ h_2, h_1}\star i)({W_{\mathcal G_1,\mathcal G_2}}^*)(\omega_{ k_2,k_1}\star i)({W_{\mathcal G_1,\mathcal G_2}}^*)Jh,k) \\
  &= (J(\omega_{ h_1, h_2}\star i)({W_{\mathcal G_1,\mathcal G_2}})^* (\omega_{ k_1,k_2}\star i)({W_{\mathcal G_1,\mathcal G_2}})^* Jh,k)\\
  &= (Jk, (\omega_{ h_1, h_2}\star i)({W_{\mathcal G_1,\mathcal G_2}})^* (\omega_{ k_1,k_2}\star i)({W_{\mathcal G_1,\mathcal G_2}})^* Jh)\\
 &=  ( (\omega_{ k_1, k_2}\star i)( W_{\mathcal G_1,\mathcal G_2} )(\omega_{ h_1,h_2}\star i)(W_{\mathcal G_1,\mathcal G_2} )Jk,Jh)  
 \end{align*}

which gives that:
 $$ (R  {_{m}\star_{s}} R)\Gamma  = \varsigma_{m,s} \Gamma \circ R $$
\end{dm}

\vskip 1cm
\subsection{The Haar operator valued weights } In this paragraph, we define two invariant operator valued weights on $ L^\infty(\mathcal G_2,\mu_2) \rtimes_{\mathfrak a} \mathfrak G(\mathcal G_1) $.

\vskip 1cm

\subsubsection{\bf Definition}
\label{enfer}
Let $T_2$     be the left    Haar operator valued weight of $\mathfrak G(\mathcal G_2)$ and  let     $T_L= \widetilde{T_2} $  be its dual operator valued weights  on $ L^\infty(\mathcal G_2,\mu_2) \rtimes_{\mathfrak a} \mathfrak G(\mathcal G_1) $  in the sense of \ref{poidsdual} and let $T_R$ be equal to $RT_LR$.

\vskip 1cm

\subsubsection{\bf Lemma}
\label{paradis}
{\it The operator valued weight  $T_L$ (resp.$T_R$) takes its  values in the range (resp.source) basis of $ L^\infty(\mathcal G_2,\mu_2) \rtimes_{\mathfrak a} \mathfrak G(\mathcal G_1)$.}
\vskip 0.2cm

\begin{dm}
Due to lemma \ref{bienfait}, $T_L =  \mathfrak a \circ T_2 \circ \mathfrak a^{-1} \circ T_{\tilde{ \mathfrak a}}= m \circ r_2^{-1} \circ T_2 \circ \mathfrak a^{-1} \circ T_{\tilde{ \mathfrak a}}$ , so  $T_L$ takes its values in $m(L^\infty(\mathcal G^0, \nu))$ which is the range basis of $ L^\infty(\mathcal G_2,\mu_2) \rtimes_{\mathfrak a} \mathfrak G(\mathcal G_1)$, due to \ref{majid2}, one has $s = R\circ m$,  so $T_R = RT_LR $ takes its values in   $s(L^\infty(\mathcal G^0, \nu))$ which gives the second part of the lemma.
\end{dm}

Hence $T_L$   is an operator valued weight from  $L^\infty(\mathcal G_2,\mu_2) \rtimes_{\mathfrak a} \mathfrak G(\mathcal G_1)$ to its range   basis, and  $\phi_L = \nu \circ m^{-1} \circ T_L$    is a dual weight, in the sense of \cite{E2} 13.1.

\vskip 1cm
\subsubsection{\bf Proposition}
\label{antipode}
{\it For any $f \in \mathcal K(\mathcal G_2), h \in \mathcal K (\mathcal G_1)$ and $\mu$-almost any $y \in \mathcal G$, one has:
 $$i)T_L(   \mathfrak a(f)(1_{s_2}{\otimes}_{r_1} \rho( h)))(y) = \int_{\mathcal G_2} f(x)h (s(x))d\lambda_2^{m(y)}(x)$$
 ii) With the notations of remark \ref{dublin}, for any $F \in \mathcal K(\mathcal  G_2   {_{s_2}\times_{r_1}} \mathcal G_1)$ and $\mu$-almost any $y \in \mathcal G$, one has:
$$T_L(\mathfrak R(F))(y) = \int_{\mathcal G_2} F(x,s(x))d\lambda_2^{m(y)}(x) $$ 
 iii) Let $ \Phi_L =  \nu \circ m^{-1} \circ T_L$ be the lifted n.sf.f weight of $T_L$ , then one has:
  $$(i \surl{\nonumber_{m} \star_{s}}_{L^{\infty}(\mathcal G^0,\nu)}\Phi_L) )\Gamma(\mathfrak a(f)(1_{s_2}{\otimes}_{r_1} \rho( h))) = T_L(\mathfrak a(f)(1_{s_2}{\otimes}_{r_1} \rho( h))). $$   
  }
\vskip 0.8cm
\begin{dm}
For any $  h \in \mathcal K (\mathcal G_1)$ and $\mu$-almost any $y \in \mathcal G$, one has:
\vskip 0.3cm

 $\begin{array}{lll} T_L(\mathfrak a(f)(1 {_{s_2}\otimes}_{r_1} \rho( h)))(y)  =  
\\

\\

 = (m \circ r_2^{-1}\circ T_2\circ \mathfrak a^{-1} \circ T_{\tilde{ \mathfrak a}})(\mathfrak a(f)(1 {_{s_2}\otimes}_{r_1} \rho( h))(y) 
\\

\\
= ( r_2^{-1}\circ T_2(f\mathfrak a^{-1} (1 {_{s_2}\otimes}_{r_1} \hat T^c_{\mathcal G_1}(\rho( h))(m(y))  &\quad \mathrm{by}\, \textit{\cite{E2}, 9.6}
\\

\\
= ( r_2^{-1}\circ T_2(f\mathfrak a^{-1} (1 {_{s_2}\otimes}_{r_1} s_1(   h_{\mid \mathcal G^0}))(m(y))  & 
\\

\\
= ( r_2^{-1}\circ T_2(fs_2(  h_{\mid \mathcal G^0}))(m(y))  &\quad \mathrm{by}\, \textit{\ref{bienfait}}

\\

\\
=\displaystyle{ \int_{\mathcal G_2} fh_{\mid \mathcal G^0}(s(x))d\lambda_2^{m(y)}(x)}

=\displaystyle{ \int_{\mathcal G_2} f(x)h(s(x))d\lambda_2^{m(y)}(x)}

\end{array}$ 

which gives i). The assertion ii) is an easy consequence of i)
\vskip 0.3cm
iii) For any $\phi \in \mathcal K(\mathcal G)$, one has:

 $\begin{array}{lll}
 \omega_\phi((i \surl{\nonumber_{m} \star_{s}}_{L^{\infty}(\mathcal G^0,\nu)}\Phi_L)\Gamma(\mathfrak a(f)(1_{s_2}{\otimes}_{r_1} \rho( h))))=
\\

\\
 
=\Phi_L(( \omega_\phi \surl{\nonumber_{m} \star_{s}}_{L^{\infty}(\mathcal G^0,\nu)}i)\Gamma(\mathfrak a(f)(1_{s_2}{\otimes}_{r_1} \rho( h))) 
\\

\\
 =\Phi_L(\mathfrak R(\Psi_{f,h}))) \hskip 6cm \mathrm{by}\, \textit{\ref{couscous 2}}
\\

\\
=\int_{\mathcal G^0}\int_{\mathcal G_2}\Psi_{f,h}(x,s(x))d\lambda_2^u(x)d\nu(u) &  \mathrm{by}\, \textit{ii)}
\\

\\
=\int_{\mathcal G^0}\int_{\mathcal G_2}\int_{\mathcal G} h(s(x))f(p_2(g^{-1})x)\delta(g^{-1})\phi\overline \phi(g^{-1}) d\lambda^{r(x))}(g)d\lambda_2^u(x)d\nu(u) 
\\

\\
=\int_{\mathcal G^0}\int_{\mathcal G_2}\int_{\mathcal G} h(s(x))f(p_2(g^{-1})x)\delta(g^{-1})|\phi|^2(g^{-1}) d\lambda^u(g)d\lambda_2^u(x)d\nu(u)   &  \mathrm{as \ \  r(x) =u}
\\

\\
=\int_{\mathcal G^0}\int_{\mathcal G_2} h(s(x))f(p_2(g^{-1})x)d\lambda_2^u(x)\int_{\mathcal G }\delta(g^{-1})|\phi|^2(g^{-1}) d\lambda^u(g)d\nu(u)    
\\

\\
=\int_{\mathcal G^0}\int_{\mathcal G }\Big (\int_{\mathcal G_2} h(s(x))f(p_2(g^{-1})x)d\lambda_2^u(x)\Big )\delta(g^{-1})|\phi|^2(g^{-1}) d\lambda^u(g)d\nu(u)   
\\

\\
= \int_{\mathcal G }\Big (\int_{\mathcal G_2} h(s(p_2(g^{-1})^{-1}x))f(x)d\lambda_2^{m(g^{-1})}(x)\Big )\delta(g^{-1})|\phi|^2(g^{-1})d\mu(g)     &   
\\

\\
\hskip 2cm \mathrm {by \ \ making \ \ the \ \ change \ \ of \ \ variable :   x\mapsto p_2(g^{-1})x}  &   
\\

\\
= \int_{\mathcal G }\Big (\int_{\mathcal G_2} h(s(x))f(x)d\lambda_2^{m(g^{-1})}(x)\Big )\delta(g^{-1})\phi(g^{-1})\overline{\phi(g^{-1})}d\mu(g)      
\\

\\
=\int_{\mathcal G}\Big (\int_{\mathcal G_2} h(s(x))f(x)d\lambda_2^{m(g)}(x)\Big ) \phi(g)\overline{\phi(g )}d\mu(g)     &  

\\

\\
\hskip 3cm \mathrm {by \ \ making \ \ the \ \ change \ \ of \ \ variable :   g \mapsto g^{-1}} 
\\

\\

= \omega_\phi(T_L(\mathfrak a(f)(1_{s_2}{\otimes}_{r_1} \rho( h))))  \hskip 4cm  \mathrm{by}\, \textit{\ref{couscous 2} and i)}

\end{array}$ 

\vskip 1cm

which ends the proof.
\end{dm}
\vskip 0.5cm

Now let's prove the left  invariance of $T_L$, using technics similar to \cite{VV} 4.12.

\subsubsection{\bf Proposition}
\label{inhonourtomien}
{\it For any positive borel fonction $\Phi$ on $\mathcal  G_2   {_{s_2}\times_{r_1}} \mathcal G_1$, one has:
$$  \int_{\mathcal G_2}\int_{\mathcal G_1} \Phi(p_2(xh),h^{-1})d\lambda_1^{s(x)}(h)d\mu_2(x) =  \int_{\mathcal G_2}\int_{\mathcal G_1} \Phi(x,h)k_L(x,h)d\lambda_1^{s(x)}(h)d\mu_2(x)$$
where $k_L(x,h) =  \displaystyle{\frac{\delta}{\delta_1}(p_1(xh))}\delta_1(h)^{-1}\delta_2(p_2(xh))\delta_2(x)^{-1}$. }

 \vskip 0.2cm
 \begin{dm}

  For any $H \in \mathcal K(\mathcal  G   {_{r}\times_{r_1}} \mathcal G_1)$, due to  Fubini's theorem, proposition \ref{use} and the left invariance for $\lambda_1$ one has:
 
 \begin{align*}
&\int_{\mathcal G_2} \int_{\mathcal G_1} \int_{\mathcal G_1}  H({g_2}^{-1}g_1,k) \delta(g_1)\delta_1(g_1)^{-1}\delta_2(g_2)^{-1}d\lambda_1^{r(g_2)}(g_1) d\lambda_1^{s(g_2)}(k)d\mu_2(g_2)\\
&= \int_{\mathcal G_2} \int_{\mathcal G_1} \int_{\mathcal G_1}  H( g_2 g_1,k) \delta(g_1)\delta_1(g_1)^{-1} d\lambda_1^{s(g_2)}(g_1) d\lambda_1^{r(g_2)}(k)d\mu_2(g_2)\\
&=\int_{\mathcal G^0} \int_{\mathcal G_2} \int_{\mathcal G_1} \int_{\mathcal G_1}  H( g_2 g_1,k) \delta(g_1)\delta_1(g_1)^{-1} d\lambda_1^{s(g_2)}(g_1) d\lambda_1^{r(g_2)}(k)d\lambda_2^u(g_2)d\nu(u)\\
&=\int_{\mathcal G^0} \int_{\mathcal G_2} \int_{\mathcal G_1} \int_{\mathcal G_1}  H( g_2 g_1,k) \delta(g_1)\delta_1(g_1)^{-1} d\lambda_1^{s(g_2)}(g_1) d\lambda_1^{u}(k)d\lambda_2^u(g_2)d\nu(u)\\
&=\int_{\mathcal G^0} \int_{\mathcal G_1} \int_{\mathcal G_2} \int_{\mathcal G_1}  H( g_2 g_1,k) \delta(g_1)\delta_1(g_1)^{-1} d\lambda_1^{s(g_2)}(g_1) d\lambda_2^u(g_2)d\lambda_1^{u}(k))d\nu(u)\\
&= \int_{\mathcal G_1} \int_{\mathcal G_2} \int_{\mathcal G_1}  H( g_2 g_1,k) \delta(g_1)\delta_1(g_1)^{-1} d\lambda_1^{s(g_2)}(g_1) d\lambda_2^{r(k)}(g_2)d\mu_1(k) \\
&= \int_{\mathcal G_1} \int_{\mathcal G}  H( g,k)  d\lambda^{r(k)}(g)d\mu_1(k) = \int_{\mathcal G_1} \int_{\mathcal G}  H( kg,k)  d\lambda^{s(k)}(g)d\mu_1(k) \\
& = \int_{\mathcal G_1}\int_{\mathcal G} \delta_1(k)^{-1}H( k^{-1}g,k^{-1})  d\lambda^{r(k)}(g)d\mu_1(k) 
 \end{align*}
 
To simplify notations let's define: $h(g_1,g_2,k) = \delta_2(g_2)^{-1}\delta_1(k)^{-1}\displaystyle{\frac{\delta}{\delta_1}(g_1)}$,  then due to proposition \ref{use}, and Fubini's theorem, this gives:

  \begin{align*}
&\int_{\mathcal G_2} \int_{\mathcal G_1} \int_{\mathcal G_1}  H({g_2}^{-1}g_1,k)  \frac{\delta}{\delta_1}(g_1) \delta_2(g_2)^{-1}d\lambda_1^{r(g_2)}(g_1) d\lambda_1^{s(g_2)}(k)d\mu_2(g_2)\\
& = \int_{\mathcal G_1}\int_{\mathcal G_2} \int_{\mathcal G_1} \delta_1(k)^{-1}H( k^{-1}g_2g_1,k^{-1})  \frac{\delta}{\delta_1}(g_1) d\lambda_1^{s(g_2)}(g_1)d\lambda_2^{r(k)}(g_2)d\mu_1(k) \\
& = \int_{\mathcal G_1}\int_{\mathcal G_2} \int_{\mathcal G_1} \delta_1(k)^{-1}H( k^{-1}g_2g_1,k^{-1}) \frac{\delta}{\delta_1}(g_1) d\lambda_1^{s(g_2)}(g_1)d\lambda_2^{r(k)}(g_2)d\lambda_1^{u}(k)d\nu(u) \\
& = \int_{\mathcal G_1}\int_{\mathcal G_2} \int_{\mathcal G_1} \delta_1(k)^{-1}H( k^{-1}g_2g_1,k^{-1}) \frac{\delta}{\delta_1}(g_1) d\lambda_1^{s(g_2)}(g_1)d\lambda_2^{u}(g_2)d\lambda_1^{u}(k)d\nu(u) \\
& = \int_{\mathcal G_2}\int_{\mathcal G_1} \int_{\mathcal G_1} \delta_1(k)^{-1}H( k^{-1}g_2g_1,k^{-1}) \frac{\delta}{\delta_1}(g_1) d\lambda_1^{s(g_2)}(g_1)d\lambda_1^{u}(k)d\lambda_2^{u}(g_2)d\nu(u) \\
& = \int_{\mathcal G_2}\int_{\mathcal G_1\times \mathcal G_1} \delta_1(k)^{-1}H( k^{-1}g_2g_1,k^{-1}) \frac{\delta}{\delta_1}(g_1) d\lambda_1^{s(g_2)}(g_1)d\lambda_1^{r(g_2)}(k)d\mu_2(g_2)\\
& = \int_{\mathcal G_2}\int_{\mathcal G_1 \times \mathcal G_1}  h(g_1,g_2,k)H( k^{-1}{g_2}^{-1}g_1,k^{-1})  d\lambda_1^{r(g_2)}(g_1)d\lambda_1^{s(g_2)}(k)d\mu_2(g_2)
 \end{align*}
 In the integral relative to $g_1$, let's use the left invariance, this gives:
 \begin{align*}
&\int_{\mathcal G_1} h(g_1,g_2,k)H( k^{-1}{g_2}^{-1}g_1,k^{-1}) d\lambda_1^{r(g_2)}(g_1) = \\
&= \int_{\mathcal G_1} h(g_1,g_2,k)H(p'_2( k^{-1}{g_2}^{-1})p'_1(k^{-1}{g_2}^{-1})g_1,k^{-1}) d\lambda_1^{r(g_2)}(g_1) \\
&= \int_{\mathcal G_1} h(p'_1(k^{-1}{g_2}^{-1})^{-1}g_1,g_2,k)H(p'_2( k^{-1}{g_2}^{-1})g_1,k^{-1}) d\lambda_1^{r(p'_1(k^{-1}{g_2}^{-1}))}(g_1) 
 \end{align*}
 but $p'_2(k^{-1}{g_2}^{-1}) = p_2(g_2k)^{-1}$ and   $p'_1(k^{-1}{g_2}^{-1}) = p_1(g_2k)^{-1}$, hence:
 
 \begin{align*}
\int_{\mathcal G_1} h(g_1,g_2,k)H( k^{-1}&{g_2}^{-1}g_1,k^{-1}) d\lambda_1^{r(g_2)}(g_1) = \\
&= \int_{\mathcal G_1} h(p_1(g_2k)g_1,g_2,k)H(p_2(g_2k)^{-1}g_1,k^{-1}) d\lambda_1^{r(p_2(g_2k))}(g_1)  \end{align*}
 
  Finally:
  \begin{align*}
&\int_{\mathcal G_2} \int_{\mathcal G_1} \int_{\mathcal G_1}  H({g_2}^{-1}g_1,k)  \frac{\delta}{\delta_1}(g_1) \delta_2(g_2)^{-1}d\lambda_1^{r(g_2)}(g_1) d\lambda_1^{s(g_2)}(k)d\mu_2(g_2)=\\
 &\int_{\mathcal G_2}\int_{\mathcal G_1}\int_{ \mathcal G_1} \Theta(g_1,g_2,k) d\lambda_1^{r(p_2(g_2k))}(g_1)d\lambda_1^{s(g_2)}(k)d\mu_2(g_2)
  \end{align*}
  where
$$ \Theta(g_1,g_2,k)  =   h(p_1(g_2k)g_1,g_2,k)H(p_2(g_2k)^{-1}g_1,k^{-1}) $$
 
 This equality can be extended to all positive borel fonctions $H$. Let us take $H = ((\frac{\delta_1}{\delta} h_1)\times (\delta_2h_2))\rho^{-1}\times h_3$, where $\rho$ is the homeomorphism $(g_1,g_2) \to g_2^{-1}g_1$ from $\mathcal  G_1   {_{r_1}\times_{r_2}} \mathcal G_2$ onto its image. Let us choose $h_1$ such that $\varphi_1: u \mapsto \int_{\mathcal G_1} h_1(g_1)d\lambda^u(g_1)$ never vanishes, so one gets:
 \begin{align*}
&\int_{\mathcal G_2} \int_{\mathcal G_1} (\varphi_1 \circ r)h_2(g_2)h_3(k) d\lambda_1^{s(g_2)}(k)d\mu_2(g_2) =\\
&= \int_{\mathcal G_2}\int_{\mathcal G_1}(\varphi_1 \circ r)h_2(p_2(g_2k))\frac{\delta}{\delta_1}(p_1(g_2k))\delta_2(g_2)^{-1}\delta_2(p_2(g_2k))\delta_1(k^{-1})\lambda_1^{s(g_2)}(k)d\mu_2(g_2)
 \end{align*}
So for any   positive borel fonction $\Phi$ on $\mathcal  G_2   {_{s_2}\times_{r_1}} \mathcal G_1$, one has:
 \begin{align*}
\int_{\mathcal G_2} \int_{\mathcal G_1}  \Phi(g_2,k)&d\lambda_1^{s(g_2)}(k)d\mu_2(g_2)= \\
 &= \int_{\mathcal G_2}\int_{\mathcal G_1}\Phi(p_2(g_2k), k^{-1})k_L(g_2,k)\lambda_1^{s(g_2)}(k)d\mu_2(g_2)
 \end{align*}
Let $\gamma : \mathcal  G_2   {_{s_2}\times_{r_1}} \mathcal G_1 \to \mathcal  G_2   {_{s_2}\times_{r_1}} \mathcal G_1$ be defined by $\gamma(g_2,k) = \gamma(p_2(g_2k),k^{-1}$, as $g_2k= p_1(g_2k)p_2(g_2k)$ hence $p_2(g_2k)k^{-1} = p_1(g_2k)^{-1}g_2$ and so $p_2(p_2(g_2k)k^{-1} = g_2
$ and $\gamma$ is a symmetry. Applying  the last formula to $\Phi \circ \gamma$ leads to the proposition.
 \end{dm}

\subsubsection{\bf Lemma}
\label{laden}
{\it The map:  $\Lambda_{\Phi_L} \mathcal R(F) \to (\delta_1k_L)^{\frac{1}{2}}F$, for $F \in\mathcal K(\mathcal  G_2   {_{s_2}\times_{r_1}} \mathcal G_1)$, realizes an isomorphism between $H_{\Phi_L}$ and $L^2(\mathcal  G_2   {_{s_2}\times_{r_1}} \mathcal G_1,\mu^2_{s_2,r_1})$.}

\vskip 0.2cm
 \begin{dm} 
 For $F \in\mathcal K(\mathcal  G_2   {_{s_2}\times_{r_1}} \mathcal G_1)$, using \ref{inhonourtomien} one has:
 \begin{align*}
 (\Lambda_{\Phi_L} \mathcal R(F) ,\Lambda_{\Phi_L} \mathcal R(F) )
 &=  \Phi_L  ( \mathcal R(F^{\#}\star  F) ) 
 = \int_{\mathcal G_2}  F^{\#}\star   F(x,s(x)) d\mu_2(x) \\
 &=  \int_{\mathcal G_2} \int_{\mathcal G_1} F^{\#}(x,h)   F(p_2(xh),h^{-1})d\lambda_1^{s(x)}(h)  )d\mu_2(u) \\
 &=  \int_{\mathcal G_2} \int_{\mathcal G_1} \overline {\delta_1^{\frac{1}{2}}F}  \delta_1^{\frac{1}{2}}F(p_2(xh),h^{-1})d\lambda_1^{s(x)}(h)  d\mu_2(x)) \\
 &=  \int_{\mathcal G_2} \int_{\mathcal G_1}  \overline {\delta_1^{\frac{1}{2}}F} \delta_1^{\frac{1}{2}}F(x,h)k_L(x,h)  d\lambda_1^{s(x)}(h)  d\mu_2(x)
 \\
 &= \mu^2_{s_2,r_1}((\delta_1k_L)^{\frac{1}{2}}F \overline{(\delta_1k_L)^{\frac{1}{2}}F} )
 \end{align*} \end{dm}
 
 So one can compute the GNS construction of $\Phi_L$ on $L^2(\mathcal  G_2   {_{s_2}\times_{r_1}} \mathcal G_1,\mu^2_{s_2,r_1})$, 
 
 \subsubsection{\bf Lemma}
 \label{nasesarko}
 {\it For $\mu^2_{s_2,r_1}$-almost any $(x,h)$ in $G_2   {_{s_2}\times_{r_1}} \mathcal G_1$, one has $k_L(p_2(xh), h^{-1}) = k_L^{-1}(x,h)$.
  }

 \vskip 0.2cm
 \begin{dm} Using the notations of proposition \ref{inhonourtomien}, one has for $\mu^2_{s_2,r_1}$-almost any $(x,h)$ in $G_2   {_{s_2}\times_{r_1}} \mathcal G_1$, one has $k_L(p_2(xh), h^{-1}) = k_L(\gamma(x,h)) = \displaystyle{\frac{d\mu^2_{s_2,r_1}\circ \gamma}{d\mu^2_{s_2,r_1}}}(\gamma(x,h))$. But $\gamma$ is a symmetry, hence: $k_L(p_2(xh), h^{-1})= \displaystyle{\frac{d\mu^2_{s_2,r_1}}{d\mu^2_{s_2,r_1}\circ \gamma}}(x,h) = k_L^{-1}(x,h)$.
 \end{dm}

\vskip 0.8cm
 \subsubsection{\bf Lemma}
 \label{tomtom}
 {\it  For any $F \in \mathcal K(\mathcal  G_2   {_{s_2}\times_{r_1}} \mathcal G_1)$, one has: $ \Delta_{\Phi_L} F =   k_L^{-1}F  $.}

 \vskip 0.2cm
 \begin{dm}
For any $F,G,H \in \mathcal K(\mathcal  G_2   {_{s_2}\times_{r_1}} \mathcal G_1)$, one has: 
\begin{align*}
( \Delta_{\Phi_L}\Lambda_{\Phi_L}\mathfrak R(F),  \Lambda_{\Phi_L}&\mathfrak R(G^{\#})) 
= ( S_{\Phi_L}\Lambda_{\Phi_L}\mathfrak R(G^{\#}), S_{\Phi_L}\Lambda_{\Phi_L}\mathfrak R(F)) \\
&= (   \Lambda_{\Phi_L}\mathfrak R(G), \Lambda_{\Phi_L}\mathfrak R(F^{\#})) = \Phi_L(\mathfrak R(F)\mathfrak R(G))\\
&= \Phi_L(\mathfrak R(F\star G )) = \int_{\mathcal G_0}\int_{\mathcal G_2}F\star G (x,s(x))d\lambda^u_2(x)d\nu(u) \\
&=  \int_{\mathcal G^0}\int_{\mathcal G_2}\int_{\mathcal G_1}F(x,h)  G(p_2(xh),h^{-1})d\lambda_1^{s(x)}(h)d\lambda^u_2(x)d\nu(u)\\
&=  \int_{\mathcal G_2}\int_{\mathcal G_1}F(x,h)  G(p_2(xh),h^{-1})d\lambda_1^{s(x)}(h)d\mu_2(x) \end{align*}
Due  proposition  \ref{inhonourtomien} and lemma \ref{nasesarko}, one has :
\begin{align*}
( \Delta_{\Phi_L}\Lambda_{\Phi_L}\mathfrak R(F),  \Lambda_{\Phi_L}&\mathfrak R(G^{\#})) 
= \int_{\mathcal G_2}\int_{\mathcal G_1}F(p_2(xh),h^{-1}) k_L(x,h) G(x,h)d\lambda_1^{s(x)}(h)d\mu_2(x) \\
&= \int_{\mathcal G_2}\int_{\mathcal G_1}F(p_2(xh),h^{-1})k_L^{-1} (p_2(xh),h^{-1}) G(x,h)d\lambda_1^{s(x)}(h)d\mu_2(x) \\
&=   \int_{\mathcal G^0}\int_{\mathcal G_2}\int_{\mathcal G_1}G(x,h) ( k_L^{-1}F) (p_2(xh),h^{-1})d\lambda_1^{s(x)}(h)d\lambda^u_2(x)d\nu(u)\\
&= \int_{\mathcal G^0}\int_{\mathcal G_2}G\star (k_L^{-1}F)(x,s(x))d\lambda^u_2(x)d\nu(u) = \Phi_L(\mathfrak R(G\star d_LF)) \\
&= \Phi_L(\mathfrak R(G)\mathfrak R(k_L^{-1}F)))= \Phi_L(\mathfrak R(G^{\#})^{\#}\mathfrak R(k_L^{-1}F))) \\
&= ( \Lambda_{\Phi_L}\mathfrak R(k_L^{-1}F),  \Lambda_{\Phi_L}\mathfrak R(G^{\#}))  
 \end{align*}
the lemma follows.
  \end{dm}

\subsubsection{\bf Remark}
\label{pmf}
As $T_R = RT_LR$ there also exists  an other density $d_R$  such that $ \Delta_{\Phi_R} F  =    d_RF  $.

\subsubsection{\bf Proposition}
\label{final}
{\it For any $F \in \mathcal K(\mathcal  G_2   {_{s_2}\times_{r_1}} \mathcal G_1)$, one has:
$$ \sigma_t^{\phi_L}(\mathfrak R(F)) = \mathfrak R(\tau^{it}F)$$
where $\tau(x,h) = \tau(x,h) = \frac{\delta}{\delta_1}(p_1(xh))\delta_2(p_2(xh)x^{-1})\delta_1(h)^{-1}$, $\mu^2_{s,r}$-almost everywhere.}

\begin{dm}
 For any $F,G \in \mathcal K(\mathcal  G_2   {_{s_2}\times_{r_1}} \mathcal G_1)$ and $\mu^2_{s,r}$-almost any $(g_2,g_1) \in  G_2   {_{s_2}\times_{r_1}} \mathcal G_1$, one has:
\begin{align*}
\sigma_t^{\phi_L}(\mathfrak R(F))&G(g_2,g_1) 
=  \Delta_{\Phi_L}^{it} \mathfrak R(F) \Delta_{\Phi_L}^{-it}G(g_2,g_1) 
 =   {k_L}^{-it} (g_2,g_1) \mathfrak R(F) \Delta_{\Phi_L}^{-it}G(g_2,g_1) \\
&= {k_L}^{-it} (g_2,g_1) \int_{\mathcal G_1} F(p_2(g_2g_1), h) \delta_1(h)^{-\frac{1}{2}} {k_L}^{it} (g_2,g_1h)G(g_2, g_1h)  d\lambda^{s(g_1)}(h) \\
&=   \int_{\mathcal G_1} \displaystyle{(\frac{k_L (g_2,g_1h)}{k_L (g_2,g_1) })^{it}}F(p_2(g_2g_1), h) \delta_1(h)^{-\frac{1}{2}} G(g_2, g_1h)  d\lambda^{s(g_1)}(h) 
\end{align*}
But one has:
\begin{align*}
\frac{k_L (g_2,g_1h)}{k_L (g_2,g_1) } 
&= \displaystyle{ \frac{ \displaystyle{\frac{\delta}{\delta_1}(p_1(g_2g_1h))}\delta_1(g_1h)^{-1}\delta_2(p_2(g_2g_1h))\delta_2(g_2)^{-1}}{\displaystyle{\frac{\delta}{\delta_1}(p_1(g_1g_2))}\delta_1(g_1)^{-1}\delta_2(p_2(g_2g_1))\delta_2(g_2^{-1})}  }\\
&= \frac{\delta}{\delta_1}(p_1(g_1g_2)^{-1}(p_1(g_2g_1h))\delta_2(p_2(g_2g_1h)p_2(g_2g_1)^{-1})\delta_1(h)^{-1}
\end{align*}
Moreover, almost everywhere: 
$$p_1(g_1g_2)^{-1}p_1(g_2g_1h)= p_1(g_1g_2)^{-1}p_1(p_1(g_2g_1)p_2(g_2g_1)h)= p_1(p_2(g_2g_1)h)$$
and:
$$ p_2(g_2g_1h)p_2(g_2g_1)^{-1} = p_2(p_1(g_2g_1)p_2(g_2g_1)h)p_2(g_2g_1)^{-1} = p_2(p_2(g_2g_1)h)p_2(g_2g_1)^{-1}$$
So:
$$ \frac{k_L (g_2,g_1h)}{k_L (g_2,g_1) } =  \tau(p_2(g_2g_1),h)$$
where:
$$ \tau(x,h) = \frac{\delta}{\delta_1}(p_1(xh))\delta_2(p_2(xh)x^{-1})\delta_1(h)^{-1}$$
The lemma follows.
\end{dm}

\subsubsection{\bf Proposition}
\label{crak}
{\it  $T_L$ is  a left invariant operator valued weight on  $L^\infty(\mathcal G_2,\mu_2) \rtimes_{\mathfrak a} \mathfrak G(\mathcal G_1)$ in the sense of definition \ref{mqg} iii).}

\begin{dm}

For any $F \in \mathcal K(\mathcal  G_2   {_{s_2}\times_{r_1}} \mathcal G_1)$ and $\mu$-almost any $y \in \mathcal G$, using the fact that $\tau(x,s(x)) = 1$ gives: 
\begin{align*}
(i \surl{\nonumber_{m} \star_{s}}_{L^{\infty}(\mathcal G^0,\nu)}\Phi_L) )\Gamma(\sigma^{\Phi_L}_t(\mathfrak R(F))(y) 
&= (i \surl{\nonumber_{m} \star_{s}}_{L^{\infty}(\mathcal G^0,\nu)}\Phi_L) )\Gamma( \mathfrak R(\tau^{it}F))(y)
\\
&= T_L(\mathfrak  R(\tau^{it}F))(y) \\
&= \int_{\mathcal G_2} (\tau^{it}F)(x,s(x))d\lambda_2^{m(y)}(x) \\
&= \int_{\mathcal G_2} F(x,s(x))d\lambda_2^{m(y)}(x) \\
&= T_L(\mathfrak  R(F))(y) \\
&= (i \surl{\nonumber_{m} \star_{s}}_{L^{\infty}(\mathcal G^0,\mu)}\Phi_L) )\Gamma( \mathfrak  R(F))(y)
\end{align*}

Hence $(i \surl{\nonumber_{m} \star_{s}}_{L^{\infty}(\mathcal G^0,\nu)}\Phi_L) )\Gamma\sigma^{T_L}_t = (i \surl{\nonumber_{m} \star_{s}}_{L^{\infty}(\mathcal G^0,\nu)}\Phi_L) )\Gamma$, so using \cite{St} theorem 6.2 and proposition \ref{antipode} iii), $T_L$ is left invariant.
\end{dm}

\subsubsection{\bf Theorem}
\label{sur vie}
{\it  The family $(L^\infty(\mathcal G^0,\nu),L^\infty(\mathcal G_2,\mu_2) \rtimes_{\mathfrak a} \mathfrak G(\mathcal G_1),m,s,\Gamma, T_L, T_R, \nu)$  is a measured quantum groupoid.}
\vskip 0.3cm
\begin{dm}
Since     $T_L$ is left invariant, then $T_R =RT_LR$ is automatically right invariant and if $\Phi_R = \nu \circ s^{-1} \circ T_R$ is the lifted weight, then using \ref{tomtom} and \ref{pmf},  $\sigma^{\Phi_R}$ and $\sigma^{\Phi_L}$ commute as these are multiplication by fonctions, the theorem follows.
  \end{dm}

\subsubsection{\bf Remark}
\label{sur vie1}
Theorem \ref{sur vie} is a generalisation of the bicrossed product construction (\cite{BSV}, \cite{VV}....)
  \vskip 2cm

\section{   Two families of examples }
\vskip1cm

	In this chapter we describe two families of examples coming from case \ref{vengeance2}  and  case \ref{vengeance3}. 
\vskip 0.5cm	
	\subsection{A matched pair of groups  action on a space } Let us use the  notations of  example \ref{vengeance2}, so $\mathcal G =X\times G$  where $G$ is a group matched pair together with a \underline{right} action on a locally compact space $X$.

	Let   $p^G_1$ and $p^G_2$ be the     almost everywhere defined   functions associated with the  matched pair $G_1,G_2$ (\cite{BSV}  3.2), then for $\nu \times  dg$ almost any $(x,g) \in X \times G$, one has:
	$$ p_1(x,g) = (x,p^G_1(g)) \  \  \ p_2(x,g) = (x.p^G_1(g),p^G_2(g))  \  \  \  m(x,g) = xp^G_1(g)$$

We shall note $a_1$ the action   of  the quantum  group $L^\infty(G_1)$ on the von Neumann algebra $L^\infty(G_2)$ coming from the usual bicrossed product construction,  and   $\Gamma_1$ the usual  coproduct  of the crossed product   $L^\infty(\ G_2) \rtimes_{  a_1} L^\infty(G_1)$. So, due to Proposition \ref{couscous 2} i), for any  $h \in  \mathcal K (  G_1)$,       any $\xi  \in \mathcal K( G\times G )$ and $dg\times dg$-almost any $(g,g') \in   G\times G$:
$$   \Gamma_1(  1  \otimes  \rho( h))\xi  (g,g') =     \int_{\mathcal G_1}h(g_1) \xi(gp_1(p_2(g')g_1),g'g_1)d g_1 $$

Thanks to  \ref{R}, one easily sees that  $L^2(X\times  G_2)\surl{\nonumber_{s_2} \otimes _{r_1}}_{L^\infty(X,\nu)}L^2(X\times   G_1)$  is isomorphic to $L^2(X\times G_2 \times  G_1, \nu \times dg_2 \times dg_1)$  (and then  to  $L^2(X\times G_2) \otimes  L^2 ( G_1)$) by the application $\theta$ such that for any $f \in \mathcal K(X\times G_2 \times X \times G_1)$ and $\nu \times dg_2 \times dg_1$-almost any $(x,g_2,g_1)$ in $X\times   G_2 \times   G_1$:
$$\theta(f)(x,g_2,g_1) = f(x,g_2,x.g_2,g_1)$$

This leads to a spatial isomorphism between $L^\infty(X\times \mathcal G_2)\surl{\nonumber_{s_2} \star_{r_1}}_{L^\infty(X,\nu)}L^\infty(X\times \mathcal G_1)$ and $L^\infty(X\times G_2  \times G_1, \nu \times dg_2 \times dg_1)$ with the same formula as for $\theta$.

So the action $\mathfrak a:  L^\infty(X\times   G_2) \to L^\infty(X\times  G_2)\surl{\nonumber_{s_2} \star_{r_1}}_{L^\infty(X,\nu)}L^\infty(X\times   G_1)$, can be identified with a one to one homomorphism $L^\infty(X\times   G_2) \to L^\infty(X\times   G_2) \otimes  L^\infty (  G_1)$
\vskip 0.5cm

\subsubsection{\bf{Remark}}
\label{metivier}
By similar arguments, in that case $L^2(\mathcal G^2_{s,m},\mu^2_{s,m})$ can be identified with the space $L^2(X\times G \times G)$ using the map $\Sigma$  such that, for any$f\in  \mathcal K( \mathcal G^2_{s,m})$, one has: $\Sigma(f)(x,g,g') = f(x,g,x.gp^G_1(g')^{-1},g')$ (one can observe that  we have $s(x,g) = m(x.gp^G_1(g')^{-1},g')$).
\vskip 1.7cm

\subsubsection{\bf{Proposition}}
\label{francemusique}
{\it i) The action $\mathfrak a$ given by \ref{bienfait}, can be identified with a usual action of $G_1$ on $L^\infty(X \times   G_2)$;  if one denotes $\tilde{\mathfrak a}$ this action, for  any $(x,g_1,g_2) \in  X\times  G_1\times G_2$ such that $g_2g_1 \in G_1G_2$ and any $f \in L^ \infty(X\times  G_2)$ :  
$$\tilde{ \mathfrak a}(f)(x,g_2,g_1) = f(x.p^G_1(g_2g_1), p^G_2(g_2g_1))$$

ii) The crossed product $L^\infty(\mathcal G_2) \rtimes_{\mathfrak a} \mathfrak G(\mathcal G_1)$ is isomorphic to the usual crossed product $L^\infty(X \times  G_2) \rtimes_{\tilde{ \mathfrak a}} L^\infty(G_1)$

iii) Using the identification of remark \ref{metivier}, for almost any $(x,g,g') \in X \times G \times G$ and any $f \in L^ \infty(X\times  G_2)$, one has
$$\Gamma(\mathfrak a(f))(x,g,g') = f(x.p_1^G(g),p^G_2(g)p_2^G(g'))  $$
 moreover, for any $h \in L^\infty(X),k \in L^\infty(G_1)$, one has:
$$  \Gamma(1 _{s_2} \otimes _{r_1} \rho(h\otimes k)) =  M(h)(1\otimes \Gamma_1(1 \otimes \rho_1(k))) $$
where $M(h)$ is (the multiplication by) the function $M(h)(x,g,g') = h(x.gp_2^G(g'))$

iv) For any $f \in \mathcal K(X\times G_2)$, any $h \in \mathcal K(X\times G_1)$ and almost any $(x,g) \in X\times G$, one has:
$$T_L(\mathfrak a(f)(1 _{s_2} \otimes _{r_1} \rho(h))(x,g) = \int_{G_2} f(xp_1^G(g),g_2)h(xp_1^G(g)g_2,e)dg_2$$}

\vskip 0.6cm

\begin{dm}
i) One easily sees that $\tilde{\mathfrak a}$ is an action.
For any $h\in L^\infty(X\times  G_2)$, any function $f \in \mathcal K(X\times G_2\times X \times G_1)$  and $\nu \times dg_2 \times dg_1$ almost any $(x,g_2,g_1)$ in $X\times \mathcal G_2 \times \mathcal G_1$, one has:
\begin{align*}
\theta(\mathfrak a (h)f)(x,g_2,g_1) 
&= \mathfrak a (h)f)(x,g_2,x.g_2,g_1) \\
&= h(p_2((x,g_2)(x.g_2,g_1))f(x,g_2,x.g_2,g_1) \\
&= h(p_2(x,g_2g_1)\theta(f)(x,g_2,g_1) \\
&=  h(x.p^G_1(g_2g_1), p^G_2(g_2g_1))\theta(f)(x,g_2,g_1)\\
&= \tilde{\mathfrak a}(h)\theta(f)(x,g_2,g_1)
\end{align*}
One deduces that  $Ad(\theta) \circ \mathfrak a  = \tilde{\mathfrak a} $, which gives i).
\vskip 0.5cm
ii) The crossed product $L^\infty(\mathcal G_2) \rtimes_{\mathfrak a} \mathfrak G(\mathcal G_1)$ is generated, in the von Neumann algebra $\mathcal L(L^2(X\times  G_2)\surl{\nonumber_{s_2} \otimes _{r_1}}_{L^\infty(X,\nu)}L^2(X\times   G_1))$,  by $\mathfrak a(L^\infty(X \times G_2))$ and $1\surl{\nonumber_{s_2} \otimes _{r_1}}_{L^\infty(X,\nu)}\mathfrak  G(\mathcal G_1)'$.

$\mathfrak  G(\mathcal G_1)'$ is generated in $\mathcal L(L^2(\mathcal G_1)$ by the image of the right regular representation of $ \mathcal G_1$, as $\mathcal G_1 = X \times  G_1$,  it is  the usual crossed product of $  L^\infty (X)$ by the right action of $G_1$, so, if one denotes by $a_1$ this action and by $\rho_1$ the right regular representation of the group $G_1$,   the von Neumann algebra $ \mathfrak  G(\mathcal G_1)'$  is generated in  $\mathcal L(L^2(X\times G_1))$ by the products $a_1(\phi)(1 \otimes \rho_1(\phi_1))$, for $\phi \in  L^\infty(X)$ and $\phi_1 \in L^\infty(G_1)$. But for  $\nu \times dg_2 \times dg_1$ almost any $(x,g_2,g_1)$ in $X\times   G_2 \times   G_1$ and any $f \in \mathcal K(X\times G_2\times X \times G_1)$ one has:

\begin{align*}
\theta\big((1_{s_2} \otimes _{r_1}[a_1(\phi)&(1 _{s_2} \otimes _{r_1} \rho_1(\phi_1))])f\big)(x,g_2,g_1) =\\
&= (1_{s_2} \otimes _{r_1}[a_1(\phi)(1 _{s_2} \otimes _{r_1} \rho_1(\phi_1))])f(x,g_2,x.g_2,g_1) \\
&= \int_{G_1}(\phi \otimes \phi_1)((x.g_2).g_1,g'_1)f(x,g_2,xg_2,g_1g'_1)dg'_1\\
&= \phi(x.(g_2g_1))\int_{G_1}  \phi_1(g'_1)\theta(f)(x,g_2,g_1g'_1)dg'_1
\end{align*}
Let $k \in L^\infty(X\times G_2)$ be defined for any $(y,g_2) \in X \times G_2$ by:
$$  k(y,g_2) = \phi(y.g_2) $$
then one has:
\begin{align*}
\theta\big((1_{s_2} \otimes _{r_1}[a_1(\phi)&(1 _{s_2} \otimes _{r_1} \rho_1(\phi_1))])f\big)(x,g_2,g_1)= \\
&= k(xp_1(g_2g_1), p_2(g_2g_1))\int_{G_1}  \phi_1(g'_1)\theta(f)(x,g_2,g_1g'_1)dg'_1 \\
&= \tilde{\mathfrak a}(k)(x,g_2,g_1)\int_{G_1}  \phi_1(g'_1)\theta(f)(x,g_2,g_1g'_1)dg'_1
\end{align*}
Hence $Ad(\theta)\circ (1_{s_2} \otimes _{r_1}[a_1(\phi)(1 _{s_2} \otimes _{r_1} \rho_1(\phi_1))] ) = \tilde{\mathfrak a}(k)(1 \otimes \rho_1(\phi_1))$, as $Ad(\theta) \circ \mathfrak a  = \tilde{\mathfrak a} $, this proves that $Ad( \theta)(L^\infty(\mathcal G_2) \rtimes_{\mathfrak a} \mathfrak G(\mathcal G_1))$  is included in $L^\infty(X \times  G_2) \rtimes_{\tilde{ \mathfrak a}} L^\infty(G_1)$ and contains $\tilde{ \mathfrak a}(L^\infty(X \times  G_2))$ and also $1 \otimes \rho_1(L^\infty(G_1))$ (using $\phi =1$), so  $Ad( \theta)$ realizes a spatial isomorphism between $L^\infty(\mathcal G_2) \rtimes_{\mathfrak a} \mathfrak G(\mathcal G_1)$ and $L^\infty(X \times  G_2) \rtimes_{\tilde{ \mathfrak a}} L^\infty(G_1)$.

iii)  Due to proposition \ref{couscous 2}, for almost any $(x,g,g') \in X \times G \times G$ and any $f \in L^ \infty(X\times  G_2)$, one has:
\begin{align*}
\Gamma(\mathfrak a(f))((x,g),&(xgp^G_1(g')^{-1},g'))   
= f(p_2(x,g)p_2(xgp^G_1(g')^{-1},g'))) \\
&= f((x.p^G_1(g),p^G_2(g)) ((x.gp^G_1(g')^{-1}).p^G_1(g'),p^G_2(g')) )\\
&= f((x.p^G_1(g),p^G_2(g)) (x.g,p^G_2(g')) )\\
&= f(x.p^G_1(g),p^G_2(g)p^G_2(g'))  
\end{align*}
Moreover, for any $h \in L^\infty(X),k \in L^\infty(G_1)$, for almost any $(x,g,g') \in X \times G \times G$, as $s(xg^Gp_1(g')^{-1},g') = x.gp_2^G(g')$ any $\xi \in L^\infty(X),  \xi' \in L^\infty(G\times G)$, as one has:

\begin{align*}
 &\Sigma \Gamma(1 _{s_2} \otimes _{r_1} \rho(h\otimes k))\Sigma^*\xi(x,g,g')  =\\
 &\ \Gamma(1 _{s_2} \otimes _{r_1} \rho(h\otimes k))\Sigma^*\xi((x,g),(xgp^G_1(g')^{-1},g'))  =\\
&= \int_{ G_1}  h(z)k(g_1)\Sigma^*\xi((x,g)p_1(p_2(xgp^G_1(g')^{-1},g')(z,g_1)),(xgp^G_1(g')^{-1},g')(z,g_1))\\& \hskip 12cm dg_1
 \end{align*}
in which $z = x.gp_2^G(g')$, hence one has:
\begin{align*}
 &\Sigma \Gamma(1 _{s_2} \otimes _{r_1} \rho(h\otimes k))\Sigma^*\xi(x,g,g')  =\\
&= \int_{ G_1}  h(z)k(g_1)\Sigma^*\xi((x,g)p_1( (xg ,p^G_2(g'))(z,g_1)),(xgp^G_1(g')^{-1},g')(z,g_1))dg_1\\
&= \int_{ G_1}  h(z)k(g_1)\Sigma^*\xi((x,g)p_1 (xg,p^G_2(g')g_1),(xgp^G_1(g')^{-1},g'g_1))dg_1\\
&= \int_{ G_1}  h(z)k(g_1)\Sigma^*\xi((x,g)  (xg,p_1^G(p^G_2(g')g_1)),(xgp^G_1(g')^{-1},g'g_1))dg_1\\
&= \int_{ G_1}  h(z)k(g_1)\Sigma^*\xi((x,p_1^G(p^G_2(g')g_1)),(xgp^G_1(g')^{-1},g'g_1))dg_1\\
&= \int_{ G_1}  h(z)k(g_1) \xi(x,p_1^G(p^G_2(g')g_1),g'g_1))dg_1\\
&= h(x.gp_2^G(g'))\int_{ G_1}  k(g_1) \xi(x,p_1^G(p^G_2(g')g_1),g'g_1))dg_1\\
&= h(x.gp_2^G(g'))(1 \otimes \Gamma_1(1\otimes \rho_1(k)))\xi(x,g,g')\\
&=M(h)(1 \otimes \Gamma_1(1\otimes \rho_1(k)))\xi(x,g,g')
\end{align*}

Where $M(h)$ is the multiplication operation by  the function:
$$M(h)(x,g,g') = h(x.gp_2^G(g'))$$
iv) Using proposition \ref{antipode} for almost any $(x,g) \in X \times G$, any $f\in \mathcal K(X\times G_1)$ and any $h \in \mathcal K(X \times G_2)$, one has:
\begin{align*}
T_L(\mathfrak a(f)(1 _{s_2} \otimes _{r_1} \rho(h))(x,g) 
&= \int_{X\times G_2} f(y,g_2)h(s(x,g_2))d\lambda_2^{m(x,g)}(y,g_2)\\
&= \int_{  G_2} f(xp_1^G(g),g_2)h(s(xp_1^G(g),g_2))dg_2\\
&= \int_{  G_2} f(xp_1^G(g),g_2)h(xp_1^G(g)g_2,e)dg_2
\end{align*}\end{dm}
\subsubsection{\bf{Remark}}
\label{hotel}
As any usual crossed product, $L^\infty(X \times  G_2) \rtimes_{\tilde{ \mathfrak a}} L^\infty(G_1)$ is isomorphic to $\widehat {\mathfrak  G}((X\times G_2)\times G_1)$ but  the   measured quantum groupoid structure we obtain for this  crossed product in \ref{francemusique} is not isomorphic to the natural  one of $\mathfrak  G((X\times G_2)\times G_1)$, recalled in \ref{R}, because for one  the basis is  $L^\infty(X )$ and for the other the basis is $L^\infty(X\times  G_2)$.

 \vskip 2cm

\subsection{The case of  a principal and transitive groupoid } 
Let's use notations similar to \ref{vengeance3},  so we suppose that $\mathcal G$ is a transitive and principal groupoid, hence of  the form: $X_1 \times X_2 \times X_1 \times X_2$ where the $X_i$'s are Hausdorff locally compact,  $ \mathcal G_1 = \underset{x_2 \in X_2}\sqcup   \  \ \  X_1\times  \{x_2\}\times X_1 \times  \{x_2\}$  and $ \mathcal G_2= \underset{x_1 \in X_1}\sqcup   \  \  \{x_1\} \times  X_2\times \{x_1\}   \times  X_2 $.  For any $(x_1,x_2,y_1,y_2)$ in $\mathcal G$, one easily sees that:
$$ p_1(x_1,x_2,y_1,y_2)= (x_1,x_2,y_1,x_2)  \ \ \ , \  \  p_2(x_1,x_2,y_1,y_2)= (y_1,x_2,y_1,y_2)$$
$$m(x_1,x_2,y_1,y_2) = (y_1,x_2)$$
\vskip0.3cm
One can identify $ \mathcal G_1$ (resp.$\mathcal G_2$) with $X_1\times X_1\times X_2$ (resp.$X_2\times X_2\times X_1$), using the map $(x_1,x_2,y_1,x_2) \mapsto (x_1, y_1,x_2)$ (resp.$(x_1,x_2,x_1,y_2) \mapsto (x_2, y_2,x_1)$); due to lemma \ref{ced}, the  Haar system of $\mathcal G_1$ is $(\delta_{x_1} \otimes \nu_1 \otimes \delta_{x_2})_{(x_1,x_2)}$. 

So   $  L^2(\mathcal G_1, \mu_1)$ (resp.$  L^2(\mathcal G_2, \mu_2)$) can be identified with $L^2(X_1\times X_1\times X_2, \nu_1\times  \nu_1\times  \nu_2)$ (resp. $L^2(X_2\times X_2\times X_1, \nu_2\times  \nu_2\times  \nu_1)$). 

This gives a spatial isomorphism between   $  L^\infty(\mathcal G_1, \mu_1)$ (resp.$  L^\infty(\mathcal G_2, \mu_2)$)  and  $L^\infty(X_1\times X_1\times X_2, \nu_1\times  \nu_1\times  \nu_2)$  (resp. $L^\infty(X_2\times X_2\times X_1, \nu_2\times  \nu_2\times  \nu_1)$). 

\vskip 0.3cm

Lemma \ref{middle1} gives  an obvious isomorphism, between the two von Neumann algebras $  L^\infty(\mathcal G_2,\mu_2)\underset{L^\infty(X_1\times X_2,\nu_1\times \nu_2)}{_{s_2}\otimes_{r_2}}L^\infty(\mathcal G_1,\mu_1)$ and  $ L^\infty(\mathcal G, \nu_1 \times \nu_2 \times\nu_1 \times\nu_2)$,  coming from the   map: $\big((x_2,y_2,x_1), (x_1,y_1,y_2)\big ) \mapsto (x_1,x_2,y_1,y_2)$. 

\vskip 0.2cm
Using this identification one has:
\begin{align*}
&\mathfrak a :  L^\infty(X_2\times X_2\times X_1,\nu_2\times \nu_2\times \nu_1) \to  \\
& L^\infty(X_2\times X_2\times X_1, \nu_2\times  \nu_2\times  \nu_1)\underset{L^\infty(X_1\times X_2,\nu_1\times \nu_2)}{_{s_2}\otimes_{r_1}} L^\infty(X_1\times X_1\times X_2, , \nu_1\times  \nu_1\times  \nu_2))
\end{align*}

and for any $f \in \mathcal K(\mathcal G_2)$, any $(x_1,x_2,y_1,y_2) \in \mathcal G$ one has: 

$$ \mathfrak a(f)\big((x_2,y_2,x_1), (x_1,y_1,y_2)\big ) = f(p_2(x_1,x_2,y_1,y_2)) = f(x_2,y_2,y_1)$$
\vskip 0.3cm
This formula can be interpreted just using the natural shift action  of the  groupoid $X_1\times X_1$ on the  fibred set $(X_1,id_{X_1})$ given for any elements $x_1,y_1$ in $ X_1$ by  $x_1 \underset{\mathfrak{s}_1}{\bf .}(x_1,y_1) = y_1$. To prove this let's give some definitions:

\vskip 0.2cm

\subsubsection{\bf{Definitions}}
\label{torture}
i) Let $\underset{\mathfrak S_1}{\bf .}$ be the  action of  $X_1\times X_1$ on the fibred set $(X_2\times X_2\times X_1, pr_3)$, where $pr_3$ is the projection on $X_1$,  defined  for any $(x_2,y_2,x_1) \in X_2\times X_2\times X_1$ and any $y_1 \in X_1$, by :
$$  (x_2,y_2,x_1) \underset{\mathfrak S_1}{\bf .}(x_1,y_1)=  (x_2,y_2,y_1)  $$
Let $(\tilde{\mathfrak  a},pr_3)$ be the corresponding action of  $\mathcal G(X_1 \times X_1)$ on $ L^\infty(X_2\times X_2 \times X_1,\nu_2\times \nu_2\times \nu_1)$ 
ii) Let $
\Sigma: L^2(X_2\times X_2\times X_1,\nu_2\times  \nu_2\times  \nu_1)\underset{L^\infty(X_1 \times X_2, \nu_1\times \nu_2)}{_{s_2}\otimes_{r_1}} L^2(X_1\times X_1\times X_2 , \nu_1\times  \nu_1\times  \nu_2)\to  \\
L^2(X_2\times X_2\times X_1,\nu_2\times \nu_2\times \nu_1)\underset{L^\infty(X_1,\nu_1)}{_{pr_3}\otimes_{r}}L^2(X_1\times X_1,\nu_1\times \nu_1)$
 be the isometric isomorphism  
  given for any $\phi  \in \mathcal K((X_2\times X_2\times X_1)\times (X_1\times X_1\times X_2)$ and almost  any $(x_1,x_2,y_1,y_2) \in \mathcal G$ by:
 $$\Sigma (\phi) \big((x_2,y_2,x_1), (x_1,y_1)\big ) = \phi\big((x_2,y_2,x_1)),(x_1,y_1,y_2)\big)$$
 iii) Let 
$
\Sigma': L^2(X_2\times X_2\times X_1,\nu_2 \times \nu_2\times \nu_1)\underset{L^\infty(X_1,\nu_1)}{_{pr_3}\otimes_{r}}L^2(X_1\times X_1,\nu_1\times \nu_1)\to  \\
 L^2(X_2\times X_2\times X_1 \times X_1,   \nu_2\times  \nu_2\times    \nu_1\times  \nu_1)
$
  be the isometric isomorphism   given for any $\phi  \in \mathcal K((X_2\times X_2\times X_1)\times (X_1\times X_1))$ and almost  any $(x_1,x_2,y_1,y_2) \in \mathcal G$ by:
 $$\Sigma' (\phi)(x_2,y_2,x_1,y_1)= \phi \big((x_2,y_2,x_1), (x_1,y_1)\big )$$
 
 \vskip 1.5cm

\subsubsection{\bf{Theorem}}
\label{franceculture}
{\it  Using the previous notations one has:
\vskip 0.3cm
i)   $\tilde{\mathfrak  a} =  Ad \Sigma \circ \mathfrak  a$,
\vskip 0.4cm
ii)   $\theta = Ad \Sigma'\Sigma$ realizes a spatial  von Neumann  isomorphism between $L^\infty(\mathcal G_2) \rtimes_{\mathfrak a} \mathfrak G(\mathcal G_1)$ and  $L^\infty(X_2^2,\nu_2^{\otimes 2}) \otimes \mathbb C1_{\mathcal L(L^2(X_1,\nu_1)}Ê\otimes  \mathcal L(L^2(X_1,\nu_1))$:  $L^\infty(\mathcal G_2) \rtimes_{\mathfrak a} \mathfrak G(\mathcal G_1)$  is isomorphic to  $L^\infty(X_2\times X_2,\nu_2\times  \nu_2)\otimes  \mathcal L(L^2(X_1,\nu_1))$ and to   $\mathfrak G(X_2Ê\times X_ 2)\otimes \hat {\mathfrak G}'(X_1\times X_1)$ ,

\vskip 0.4cm
iii) \hskip 0.1cm If $\tau$: $L^2(X_2^3\times X_1^3,\nu_2^{ 3}\times \nu_1^{  3}) \to L^2(X_2^2\times X_1\times X_2 \times X_1^2,\nu_2^{ 2}\times \nu_1\otimes \nu_2\times \nu_1^{ 2})$ is the map which flips the third and fourth factor, one has: 
$$(\theta _s\star_m \theta) \Gamma \theta^* = Ad \tau (\Gamma_{ X_2^2} \otimes \widehat\Gamma'_{X_1^2})$$ 
}
\vskip 0.8cm
\begin{dm}
\vskip 0.3cm
i) The first assertion  is obvious.

ii) For any $f_1,h_1,k_1 \in \mathcal K(X_1)$, any $f_2,g_2 \in \mathcal K(X_2)$, any $\xi \in \mathcal K(X_2\times X_2 \times X_1 \times X_1)$, and almost any $(x_2,y_2,x_1,y_1)\in X_2\times X_2 \times X_1 \times X_1$, due to i)  one has:
 $$\theta [ \mathfrak a(f_2\otimes g_2 \otimes f_1)]\xi(x_2,y_2,x_1,y_1) = (f_2 \otimes g_2 \otimes 1 \otimes f_1)(x_2,y_2,x_1,y_1)\xi(x_2,y_2,x_1,y_1) $$ 
 so:
  $$\theta [ \mathfrak a(f_2\otimes g_2 \otimes f_1)] =  f_2 \otimes g_2 \otimes 1 \otimes f_1$$
  Let's denote $T_{\psi_1}$, for any $\psi_1 \in \mathcal K(X_1\times X_1)$, the integral (compact) operator defined for any $\xi_1\in \mathcal K(X_1)$  and almost any $x_1$ in $X_1$  by: 
  $$(T_{\psi_1}\xi_1)(x_1) = \int_{X_1} \psi_1(x_1,z_1)\xi_1(z_1)d\nu_1(z_1)$$
\vskip 0.5cm

A straightforward calculation gives that: 
 
\begin{align}
\theta [ \mathfrak a(f_2\otimes g_2 \otimes f_1 )(1 \underset{\nu}{_{s_2}\otimes _{r_1}}\rho(h_1\otimes k_1\otimes  h_2)] =  f_2\otimes g_2h_2\otimes 1\otimes T_{f_1h_1\otimes k_1}  
\end{align}

 The assertion ii) follows.

iii) One easily sees that the coproduct $\Gamma_{  \mathcal G_2 } : L^\infty(\mathcal G_2) \to L^\infty(\mathcal G_2)_s\star_r L^\infty(\mathcal G_2)$ is given for any $f_2,g_2  \in   X_2$ and any $f_1 \in \mathcal K(X_1)$ by:
$$ \Gamma_{  \mathcal G_2 }(f_2\otimes g_2\otimes f_1) = (f_2\otimes 1\otimes f_1)_s\otimes_r (1\otimes g_2 \otimes 1)$$
Moreover, using \ref{bienfait}, one has:
\begin{align*}
&(\theta\circ m)(f_1 \otimes f_2) = ( \theta\circ \mathfrak a\circ r_2)(f_1\otimes   f_2))) =Ê\theta(\mathfrak a(f_2 \otimes 1 \otimes f_1) = f_2 \otimes 1 \otimes 1 \otimes f_1
\\
 &(\theta\circ s)(f_1 \otimes f_2)  = ( \theta\circ \mathfrak a\circ s_2)(f_1\otimes   f_2) = \theta(\mathfrak a(1 \otimes f_2 \otimes f_1) = 1 \otimes f_2 \otimes 1 \otimes f_1
  \end{align*} 

This obviously gives an isometric isomorphism of  the hilbert space  $L^2(X_2^2\times X_1^2,\nu_2^{\otimes 2}\otimes \nu_1^{\otimes 2}) \underset{L^\infty(X_1\times X_2,\nu_1\otimes \nu_2)}{_{\theta \circ s}\otimes_{\theta \circ m}} L^2(X_2^2\times X_1^2,\nu_2^{\otimes 2}\otimes \nu_1^{\otimes 2})$ onto the hilbert space $L^2(X_2^2\times X_1\times X_2 \times X_1^2,\nu_2^{\otimes 2}\otimes \nu_1\otimes \nu_2\otimes \nu_1^{\otimes 2})$, if $\Psi$ is this map,  for any $\xi^1,  \xi^2  \in \mathcal K(X_2^2\times X_1^2)$ one has:  
\begin{align}
 \Psi(\xi^1{_{\theta \circ s}\otimes_{\theta \circ m}}\xi^2)(x_2,y_2,x_1,z_2,y_1,z_1) = \xi^1(x_2,y_2,x_1,y_1)\xi^2(y_2,z_2,z_1,y_1)
 \end{align}
Hence, using \ref{couscous1}:
\begin{align*}
(\theta _s\star_m \theta) \Gamma \theta^*(f_2\otimes g_2\otimes 1 \otimes 1)
&=  (f_2\otimes 1\otimes 1 \otimes 1)_{\theta \circ s}\otimes_{\theta \circ m}  (1\otimes g_2 \otimes 1\otimes 1) 
\end{align*}
due to  (2) and \ref{fima1}, this gives:
\begin{align*}
(\theta _s\star_m \theta) \Gamma \theta^*(f_2\otimes g_2\otimes 1 \otimes 1)
&=Ad \tau (\Gamma_{ X_2^2 } \otimes \widehat\Gamma'_{ X_1^2 })(f_2\otimes g_2\otimes 1 \otimes 1)
\end{align*}
Quite simple computations, also imply that for any $h_1,k_1 \in \mathcal K(X_1)$:
\begin{align*}
 (\theta _s\star_m \theta) \Gamma& (1 \underset{\nu}{_{s_2}\otimes _{r_1}}\rho(h_1\otimes k_1\otimes  1))   =   Ad \tau (\Gamma_{ X_2^2 } \otimes \widehat\Gamma'_{ X_1^2 })  \theta^*(1 \underset{\nu}{_{s_2}\otimes _{r_1}}\rho(h_1\otimes k_1\otimes  1)) 
\end{align*}
which gives iii). 

iv) Due to (1), \ref{antipode} and \ref{fima1}, for any $f_1,g_1\in \mathcal K(X_1)$ any $f_2,g_2 \in \mathcal K(X_2)$,  one has:
\begin{align*}
\theta T_L\theta^*&(  f_2\otimes g_2\otimes 1 \otimes  T_{h_1\otimes k_1})     = (T_{X_2^2} \otimes \widehat T_{X_1^2}')(f_2\otimes g_2 \otimes 1 \otimes T_{h_1\otimes k_1})
\end{align*}

Which gives : $\theta T_L\theta^* =  T_{X_2^2} \otimes \widehat T_{X_1^2}' $.\end{dm}

\end{document}